\documentclass[12pt]{article}

\usepackage{amsmath}
\usepackage{amsthm}
\usepackage{amssymb}
\usepackage[left=2.5cm,right=2.4cm,top=1.5cm,bottom=2.0cm]{geometry}
\usepackage{tikz}
\usetikzlibrary{positioning,shapes,calc,backgrounds,fit,matrix}
\usepackage{pgfplots}
\usetikzlibrary{plotmarks}
\usepackage{algpseudocode}
\usepackage{algorithm}
\usepackage{hyperref}

\theoremstyle{definition}
\newtheorem{Rem}{Remark}

\def\eps{\varepsilon}

\newcommand*\sq{\mathbin{\vcenter{\hbox{\rule{1ex}{1ex}}}}}

\usepackage{tikz} \usetikzlibrary{positioning,shapes,calc,backgrounds,fit,plotmarks}
\usepackage{pgfplots}
\pgfplotsset{compat=newest,
every axis/.append style={axis x line=bottom,
                          scale only axis,
                          },
}

\date{June 22, 2018}

\begin{document}
\title{Preconditioners and Tensor Product Solvers for Optimal Control Problems from Chemotaxis}
\author{Sergey Dolgov\footnote{University of Bath, Claverton Down, BA2 7AY, Bath, United Kingdom. {\tt s.dolgov@bath.ac.uk}} and John W. Pearson\footnote{School of Mathematics, The University of Edinburgh, James Clerk Maxwell Building, The King's Buildings, Peter Guthrie Tait Road, Edinburgh, EH9 3FD, United Kingdom. {\tt j.pearson@ed.ac.uk}}}

\maketitle

\begin{abstract}
In this paper, we consider the fast numerical solution of an optimal control formulation of the Keller--Segel model for bacterial chemotaxis. Upon discretization, this problem requires the solution of huge-scale saddle point systems to guarantee accurate solutions. We consider the derivation of effective preconditioners for these matrix systems, which may be embedded within suitable iterative methods to accelerate their convergence. We also construct low-rank tensor-train techniques which enable us to present efficient and feasible algorithms for problems that are finely discretized in the space and time variables. Numerical results demonstrate that the number of preconditioned GMRES iterations depends mildly on the model parameters.
Moreover, the low-rank solver makes the computing time and memory costs sublinear in the original problem size.
\end{abstract}


\section{Introduction}\label{sec:Intro}

The process of chemotaxis in biology describes the movement of cells or organisms in a directed fashion as a response to external chemical signals. In 1971, Keller and Segel presented a mathematical model for bacterial chemotaxis \cite{KellerSegel}. In essence, for large numbers of bacteria, it is predicted that the bacteria will on average move up gradients of the chemoattractant concentration.

Since Keller and Segel's work, an area of numerical mathematics that has become a subject of significant interest is that of PDE-constrained optimization, where one wishes to predict the circumstances in which some physical (or in this case biological) objective occurs, subject to a system of PDEs describing the process. Using this technology, one is able to pose an inverse problem for the chemotaxis mechanism: given an observed bacterial cell concentration profile, what can be said about the external chemoattractant at the boundaries of a domain of interest? The constraints for this problem are therefore the PDEs describing bacterial chemotaxis. This is a parameter identification problem that has been considered in literature such as \cite{LBP,Potschka}, and in particular it was shown numerically by Lebiedz and Brandt-Pollmann that ``it is possible to systematically control spatiotemporal dynamical behavior'' \cite{LBP}.

The fast and efficient iterative solution of PDE-constrained optimization problems has increasingly become an active area of research, and in particular it is now widely recognised that the incorporation of effective preconditioners to accelerate iterative schemes is highly beneficial from a computational point-of-view. Preconditioning theory and numerics for a number of steady \cite{PSW_StateConstraints,PW2,PW1,RDW,SchoberlZulehner,Zulehner} and time-dependent \cite{DPSS,PearsonStoll,PSW,SPM} problems have been established, with \cite{PearsonStoll,SPM} describing the resulting solvers for reaction--diffusion problems from chemistry and biology. In this paper, we derive a potent preconditioner for the chemotaxis problem based on the saddle point structure of the matrix systems resulting from Newton-type iterations of the nonlinear PDEs.

When solving these optimization problems, which often involve the solution of a system of PDEs with initial conditions coupled with adjoint PDEs equipped with final-time conditions, there are many challenges arising from the time-dependent component of the problem in particular, due to the forward-backward solves required, and the associated scaling of computational complexity with the fineness of the grid in the time variable.
Difficulties also arise from nonlinear problems, due to the matrices arising from the PDE system varying in structure at every time step, unlike linear problems for which some matrices can be re-used repeatedly within a solver.
For time-dependent nonlinear problems that arise from chemotaxis, computer storage is therefore a significant bottleneck, unless a numerical algorithm is specifically tailored in order to mitigate this.

To combat this issue, in addition to presenting our new preconditioner, we describe an approach for approximating the solution of our problem in a low-rank format, namely the Tensor Train decomposition \cite{osel-tt-2011}.
Low-rank tensor techniques emerge from the separation of variables and the Fourier method for solving PDEs.
We can approximate the solution in the form $z(x,y) \approx \sum_{\alpha=1}^{r} v_{\alpha}(x) w_{\alpha}(y)$, using a possibly small number of terms $r$.
In this case, the discretized univariate functions in the low-rank decomposition are much cheaper to store than the original multivariate function.
The discretized separation of variables requires the low-rank approximation of matrices (for two variables $x,y$), or tensors (for three or more variables).
Practical low-rank tensor algorithms employ robust tools of linear algebra, such as the singular value decomposition, to deliver an optimal low-rank approximation for a desired accuracy.
Extensive reviews on the topic can be found in \cite{hackbusch-2012,bokh-surv-2015}.

The efficiency of low-rank decompositions depends crucially on the value of the rank $r$, which in turn reflects the structure of a function.
Discontinuous functions, in particular level set functions, may require high ranks.
However, smooth functions allow very accurate low-rank approximations, and hence a sublinear complexity of the inverse problem solution \cite{uschmajew-approx-rate-2013,tee-tensor-2003}.
The inverse problem implies driving the solution to a desired state,
which usually has a simple (and hence low-rank) structure.
Therefore, as long as we avoid discontinuous functions in our formulation,
the low-rank techniques can be very efficient for the inverse problem.
This is demonstrated in our computational experiments.

This paper is structured as follows.
In Section \ref{sec:Problem} we describe the problem statement of which the numerical solution is considered.
In Section \ref{sec:Matrix} we present the structure of the matrix systems that result from the discretization of the system of PDEs.
In Section \ref{sec:Preconditioner} we present our preconditioning strategy for these systems, with numerical results provided in Section \ref{sec:NumEx1}.
We describe the low-rank tensor decomposition which is employed for the matrix systems in Section \ref{sec:LowRank}, with additional numerical experiments relating to this approach carried out in Section \ref{sec:NumEx2}.
Finally, concluding remarks are made in Section \ref{sec:Conc}.

\section{Problem statement}\label{sec:Problem}

We examine the following problem describing the optimal control of a bacterial chemotaxis system, based on studies in literature such as \cite{LBP} and \cite[Chapter 13]{Potschka}:
\begin{align}
\ \label{CostFunctional} \min_{z,c,u}~~\frac{1}{2}\int_{\Omega}\left(z(\mathbf{x},T)-\widehat{z}\right)^2+\frac{\gamma_c}{2}\int_{\Omega}\left(c(\mathbf{x},T)-\widehat{c}\right)^2+\frac{\gamma_u}{2}{}&{}\int_{\partial\Omega\times(0,T)}u^2
\end{align}
subject to
\begin{align*}
\frac{\partial{}z}{\partial{}t}-D_{z}\nabla^{2}z-\alpha\nabla\cdot\left(\frac{\nabla{}c}{(1+c)^2}z\right)=0\quad\quad&\text{on }\Omega\times(0,T), \\
\ \nonumber \frac{\partial{}c}{\partial{}t}-\nabla^{2}c+\rho{}c-w\frac{z^2}{1+z^2}=0\quad\quad&\text{on }\Omega\times(0,T), \\
\end{align*}
equipped with the boundary conditions and initial conditions:
\begin{align}
\ \nonumber \frac{\partial{}z}{\partial{}n}=0\hspace{1.8em}\quad&\text{on }\partial\Omega\times(0,T), \\
\ \nonumber \frac{\partial{}c}{\partial{}n}+\beta{}c=\beta{}u\hspace{1.1em}\quad&\text{on }\partial\Omega\times(0,T), \\
\ \nonumber z(\mathbf{x},0)=z_0(\mathbf{x})\quad&\text{on }\Omega, \\
\ \nonumber c(\mathbf{x},0)=c_0(\mathbf{x})\quad&\text{on }\Omega.
\end{align}
This problem is solved on a space-time domain $\Omega\times(0,T)$ with boundary $\partial\Omega\times(0,T)$, and for $\Omega\subset\mathbb{R}^2$. The variables $z$, $c$ denote \emph{state variables}, corresponding to the bacterial cell density and chemoattractant concentration respectively, with $u$ the \emph{control variable}, $\widehat{z}$, $\widehat{c}$ given \emph{desired states}, $z_0$, $c_0$ given initial conditions, and $\gamma_c$, $\gamma_u$, $D_z$, $\alpha$, $\rho$, $w$, $\beta$ given (positive) parameters. We highlight that, by construction of the problem, the control $u$ in some sense relates to the gradient of chemoattractant concentration on the boundary of the domain of interest. The form of the boundary condition which enforces the control makes this a \emph{boundary control problem}. In this PDE-constrained optimization model, we wish to discover what the profile of this control must be in order for the biological system to behave in a way prescribed by the desired states $\widehat{z}$, $\widehat{c}$.

\begin{Rem}
Although derived similarly, the main difference between the works of \cite{LBP} and \cite{Potschka} is that \cite{LBP} considers solely the misfit between $z$ and $\widehat{z}$, regularized by a term involving the final time $T$. We believe that the methods introduced in this paper are equally applicable to either cost functional.
\end{Rem}

We now consider the first and second derivatives of the Lagrangian\footnote{For ease of notation, we exclude initial conditions within the definition of the Lagrangian.}
\begin{align*}
\ \mathcal{L}(z,c,u,p,q)={}&\frac{1}{2}\int_{\Omega}\left(z(\mathbf{x},T)-\widehat{z}\right)^2+\frac{\gamma_c}{2}\int_{\Omega}\left(c(\mathbf{x},T)-\widehat{c}\right)^2+\frac{\gamma_u}{2}\int_{\partial\Omega\times(0,T)}u^2 \\
\ &\quad\quad+\int_{\Omega\times(0,T)}p_{\Omega}\left(\frac{\partial{}z}{\partial{}t}-D_{z}\nabla^{2}z-\alpha\nabla\cdot\left(\frac{\nabla{}c}{(1+c)^2}z\right)\right) \\
\ &\quad\quad+\int_{\Omega\times(0,T)}q_{\Omega}\left(\frac{\partial{}c}{\partial{}t}-\nabla^{2}c+\rho{}c-w\frac{z^2}{1+z^2}\right) \\
\ &\quad\quad+\int_{\partial\Omega\times(0,T)}p_{\partial\Omega}\left(\frac{\partial{}z}{\partial{}n}\right)+\int_{\partial\Omega\times(0,T)}q_{\partial\Omega}\left(\frac{\partial{}c}{\partial{}n}+\beta{}c-\beta{}u\right),
\end{align*}
where $p$ and $q$ denote the adjoint variables corresponding to $z$ and $c$, with $p_{\Omega}$, $q_{\Omega}$ the components of $p$, $q$ within the interior of $\Omega$, and $p_{\partial\Omega}$, $q_{\partial\Omega}$ the components on the boundary. We arrive at the following system for the Newton formulations of the first-order optimality conditions:
\begin{align}
\ \label{Newton1} &\frac{\partial{}s_z}{\partial{}t}-D_{z}\nabla^{2}s_z+\alpha\nabla\cdot\left(\nabla\left(\frac{1}{1+c}\right)s_z\right)+\alpha\nabla\cdot\left(\nabla\left(\frac{1}{(1+c)^2}s_c\right)z\right) \\
\ \nonumber &\quad\quad\quad\quad=-\left(\frac{\partial{}z}{\partial{}t}-D_{z}\nabla^{2}z-\alpha\nabla\cdot\left(\frac{\nabla{}c}{(1+c)^2}z\right)\right)\quad\text{on }\Omega\times(0,T), \\
\ \label{Newton2} &\frac{\partial{}s_c}{\partial{}t}-\nabla^{2}s_c+\rho{}s_c-2w\frac{z}{(1+z^2)^2}s_z \\
\ \nonumber &\quad\quad\quad\quad=-\left(\frac{\partial{}c}{\partial{}t}-\nabla^{2}c+\rho{}c-w\frac{z^2}{1+z^2}\right)\hspace{3.45em}\quad\text{on }\Omega\times(0,T),
\\
\ \label{Newton3} &\gamma_u{}s_u-\beta{}s_q=-\left(\gamma_u{}u-\beta{}q\right)\hspace{10.7em}\quad\text{on }\partial\Omega\times(0,T), \\
\ \label{Newton4} &\chi_{\Omega_T}(s_z)-2wq\frac{1-3z^2}{(1+z^2)^3}s_z+\alpha\nabla\left(\frac{2c}{(1+c^2)^2}s_c\right)\cdot\nabla{}p \\
\ \nonumber &\quad\quad\quad\quad-\frac{\partial{}s_p}{\partial{}t}-D_{z}\nabla^{2}s_p-\alpha\nabla\left(\frac{1}{1+c}\right)\cdot\nabla{}s_p-2w\frac{z}{(1+z^2)^2}s_q \\
\ \nonumber &\quad\quad\quad\quad=\widehat{z}-\left(\chi_{\Omega_T}(z)-\frac{\partial{}p}{\partial{}t}-D_{z}\nabla^{2}p-\alpha\nabla\left(\frac{1}{1+c}\right)\cdot\nabla{}p-2w\frac{zq}{(1+z^2)^2}\right) \\
\ \nonumber &\quad\quad\quad\quad\hspace{19.55em}\text{on }\Omega\times(0,T), \\
\ \label{Newton5} &-\alpha{}p\nabla\cdot\left(\nabla\left(\frac{1}{(1+c)^2}\right)s_z\right)+\gamma_{c}\chi_{\Omega_T}(s_c)+\alpha{}p\nabla\cdot\left(\nabla\left(\frac{2c}{(1+c^2)^2}s_c\right)z\right) \\
\ \nonumber &\quad\quad\quad\quad-\alpha\nabla\cdot\left(\nabla\left(\frac{1}{(1+c)^2}\right)z\right)s_p-\frac{\partial{}s_q}{\partial{}t}-\nabla^{2}s_q+\rho{}s_q \\
\ \nonumber &\quad\quad\quad\quad=\gamma_{c}\widehat{c}-\left(\gamma_{c}\chi_{\Omega_T}(c)-\alpha{}p\nabla\cdot\left(\nabla\left(\frac{1}{(1+c)^2}\right)z\right)-\frac{\partial{}q}{\partial{}t}-\nabla^{2}q+\rho{}q\right) \\
\ \nonumber &\quad\quad\quad\quad\hspace{19.55em}\text{on }\Omega\times(0,T),
\end{align}
where $s_z$, $s_c$, $s_u$, $s_p$, $s_q$ are the Newton updates for $z$, $c$, $u$, $p$, $q$, and $\chi_{\Omega_T}(\sq)$ denotes a function that restricts the variable to time $t=T$. The boundary conditions for the state and adjoint variables are given by
\begin{eqnarray*}
\ \frac{\partial{}s_z}{\partial{}n}=0,\quad\frac{\partial{}s_c}{\partial{}n}+\beta{}s_c=\beta{}s_u,\quad\frac{\partial{}s_p}{\partial{}n}=0,\quad\frac{\partial{}s_q}{\partial{}n}+\beta{}s_q=0\quad\text{on }\partial\Omega\times(0,T),
\end{eqnarray*}
with initial and final-time conditions
\begin{eqnarray*}
\ s_z(\mathbf{x},0)=0,\quad{}s_c(\mathbf{x},0)=0,\quad{}s_p(\mathbf{x},T)=0,\quad{}s_q(\mathbf{x},T)=0\quad\text{on }\Omega,
\end{eqnarray*}
assuming an initial guess is chosen that satisfies the initial conditions for $z$, $c$ and the final-time conditions for $p$, $q$.

\begin{Rem}
We highlight that there also exist chemotaxis problems which may be written in distributed control form. For example the work in \cite{EPS}, on the identification of chemotaxis models with volume-filling, considers (amongst others) a problem which may be interpreted in our setting in the following way:
\begin{align*}
\ \min_{z,f}~~\frac{1}{2}\int_{\Omega\times(0,T)}\left(z-\widehat{z}\right)^2+\frac{\gamma}{2}\int_{\Omega\times(0,T)}\left[f^2+|\nabla{}f|^2\right]& \\
\ \emph{s.t.}\quad\frac{\partial{}z}{\partial{}t}-\nabla^{2}z+f\nabla^{2}c+\nabla{}f\cdot\nabla{}c=0\hspace{1.8em}\quad\emph{on }&\Omega\times(0,T), \\
\ -\nabla^{2}c+c=z\hspace{1.8em}\quad\emph{on }&\Omega\times(0,T), \\
\ \frac{\partial{}z}{\partial{}n}-f\frac{\partial{}c}{\partial{}n}=0\hspace{1.8em}\quad\emph{on }&\partial\Omega\times(0,T), \\
\ \frac{\partial{}c}{\partial{}n}=0\hspace{1.8em}\quad\emph{on }&\partial\Omega\times(0,T), \\
\ z(\mathbf{x},0)=z_0(\mathbf{x})\quad\emph{on }&\Omega,
\end{align*}
where $\gamma$ is a positive constant, and $f(z)$ denotes the chemoattractant sensitivity. The challenge in this case is to discover the necessary profile of the function $f$ in order to drive the chemoattractant to a particular state. We believe that variants of the techniques introduced in this paper could also be applied to this distributed control problem.
\end{Rem}

\section{Matrix systems for Newton and Gauss--Newton}\label{sec:Matrix}

In this section, we describe the matrix systems which are obtained by discretization of the optimization problem \eqref{CostFunctional} using the finite element method.

Concatenating the Newton equations \eqref{Newton1}--\eqref{Newton5}, along with boundary conditions and initial/final-time conditions, gives a block matrix system of the following form:
\begin{align}
\ \label{Newton} &\left[\begin{array}{ccccc}
\mathcal{L}_{zz} & \mathcal{L}_{zc} & 0 & \mathcal{L}_{zp} & \mathcal{L}_{zq} \\
\mathcal{L}_{cz} & \mathcal{L}_{cc} & 0 & \mathcal{L}_{cp} & \mathcal{L}_{cq} \\
0 & 0 & \gamma_{u}\cdot\text{Id} & 0 & -\beta\chi_{\partial\Omega}(\sq)^\top \\
\mathcal{L}_{pz} & \mathcal{L}_{pc} & 0 & 0 & 0 \\
\mathcal{L}_{qz} & \mathcal{L}_{qc} & -\beta\chi_{\partial\Omega}(\sq) & 0 & 0 \\
\end{array}\right]\left[\begin{array}{c}
s_z \\
s_c \\
s_u \\
s_p \\
s_q \\
\end{array}\right] \\
\ \nonumber &\quad\quad\quad\quad=\left[\begin{array}{c}
\widehat{z}-\left(\chi_{\Omega_T}(z)-\frac{\partial{}p}{\partial{}t}-D_{z}\nabla^{2}p-\alpha\nabla\left(\frac{1}{1+c}\right)\cdot\nabla{}p-2w\frac{zq}{(1+z^2)^2}\right) \\
\gamma_{c}\widehat{c}-\left(\gamma_{c}\chi_{\Omega_T}(c)-\alpha{}p\nabla\cdot\left(\nabla\left(\frac{1}{(1+c)^2}\right)z\right)-\frac{\partial{}q}{\partial{}t}-\nabla^{2}q+\rho{}q\right) \\
-\left(\gamma_u{}u-\beta{}q\right) \\
-\left(\frac{\partial{}z}{\partial{}t}-D_{z}\nabla^{2}z-\alpha\nabla\cdot\left(\frac{\nabla{}c}{(1+c)^2}z\right)\right) \\
-\left(\frac{\partial{}c}{\partial{}t}-\nabla^{2}c+\rho{}c-w\frac{z^2}{1+z^2}\right) \\
\end{array}\right],
\end{align}
where
\begin{align*}
\ \left[\begin{array}{cc}
\mathcal{L}_{zz} & \mathcal{L}_{zc} \\
\mathcal{L}_{cz} & \mathcal{L}_{cc} \\
\end{array}\right]={}&\left[\begin{array}{cc}
\chi_{\Omega_T}(\sq)-2wq\frac{1-3z^2}{(1+z^2)^3} & \alpha\nabla\left(\frac{2c}{(1+c^2)^2}{\sq}\right)\cdot\nabla{}p \\
-\alpha{}p\nabla\cdot\left(\nabla\left(\frac{1}{(1+c)^2}\right)\sq\right) & \gamma_{c}\chi_{\Omega_T}(\sq)+\alpha{}p\nabla\cdot\left(\nabla\left(\frac{2c}{(1+c^2)^2}\sq\right)z\right) \\
\end{array}\right], \\
\ \left[\begin{array}{cc}
\mathcal{L}_{pz} & \mathcal{L}_{pc} \\
\mathcal{L}_{qz} & \mathcal{L}_{qc} \\
\end{array}\right]={}&\left[\begin{array}{cc}
\frac{\partial}{\partial{}t}-D_{z}\nabla^{2}+\alpha\nabla\cdot\left(\nabla\left(\frac{1}{1+c}\right)\sq\right) & \alpha\nabla\cdot\left(\nabla\left(\frac{1}{(1+c)^2}\sq\right)z\right) \\
-2w\frac{z}{(1+z^2)^2} & \frac{\partial}{\partial{}t}-\nabla^{2}+\rho\cdot\text{Id} \\
\end{array}\right], \\
\ \left[\begin{array}{cc}
\mathcal{L}_{zp} & \mathcal{L}_{zq} \\
\mathcal{L}_{cp} & \mathcal{L}_{cq} \\
\end{array}\right]={}&\left[\begin{array}{cc}
-\frac{\partial}{\partial{}t}-D_{z}\nabla^{2}-\alpha\nabla\left(\frac{1}{1+c}\right)\cdot\nabla &-2w\frac{z}{(1+z^2)^2} \\
-\alpha\nabla\cdot\left(\nabla\left(\frac{1}{(1+c)^2}\right)z\right) & -\frac{\partial}{\partial{}t}-\nabla^{2}+\rho\cdot\text{Id} \\
\end{array}\right],
\end{align*}
with $\text{Id}$ denoting the identity operator, and $\chi_{\partial\Omega}(\sq)$ representing a function restricted to the boundary $\partial\Omega$.

As an alternative to solving the Newton system \eqref{Newton}, it is possible to instead consider a Gauss--Newton approximation, where one neglects second derivatives within the $(1,1)$-block of the saddle point matrix as defined in Section \ref{sec:Preconditioner}. This results in the solution of systems
\begin{eqnarray}
\ \label{GN} \left[\begin{array}{ccccc}
\chi_{\Omega_T}(\sq) & 0 & 0 & \mathcal{L}_{zp} & \mathcal{L}_{zq} \\
0 & \gamma_{c}\chi_{\Omega_T}(\sq) & 0 & \mathcal{L}_{cp} & \mathcal{L}_{cq} \\
0 & 0 & \gamma_{u}\cdot\text{Id} & 0 & -\beta\chi_{\partial\Omega}(\sq)^\top \\
\mathcal{L}_{pz} & \mathcal{L}_{qz} & 0 & 0 & 0 \\
\mathcal{L}_{pc} & \mathcal{L}_{qc} & -\beta\chi_{\partial\Omega}(\sq) & 0 & 0 \\
\end{array}\right]\left[\begin{array}{c}
s_z \\
s_c \\
s_u \\
s_p \\
s_q \\
\end{array}\right]=\mathbf{b},
\end{eqnarray}
where $\mathbf{b}$ is the same right-hand side vector as in \eqref{Newton}.

To be more explicit about the $\chi_{\Omega_T}(\sq)$ and $\chi_{\partial\Omega}(\sq)$ terms, the associated matrices contain entries of the form $\int_{\Omega}\phi_i\cdot\phi_j|_{t=T}$ and $\int_{\partial\Omega\times(0,T)}\phi_i\cdot\phi_j|_{\partial\Omega}$ respectively, for finite element basis functions $\{\phi_i\}$ of the same form for each PDE variable.

\subsection{Additional control constraints}\label{sec:Matrix_ControlConstraints}

It is perfectly reasonable to add the following control constraint:
\begin{eqnarray*}
\ u_-(\mathbf{x},t)\leq{}u\leq{}u_+(\mathbf{x},t)\quad\text{a.e. on }\partial\Omega\times(0,T),
\end{eqnarray*}
for given functions $u_-$, $u_+$, into the PDE-constrained optimization model. In other words, we prescribe that the chemoattractant must behave in a ``sensible'' (physical) way on the boundary of the domain of interest. One way in which we can tackle this additional term is to modify the cost functional \eqref{CostFunctional} to add a Moreau--Yosida regularization term (see \cite{ItoKunisch}) for the bound constraints, thereby minimizing instead
\begin{align*}
\ \min_{z,c,u}~~&\frac{1}{2}\int_{\Omega}\left(z(\mathbf{x},T)-\widehat{z}\right)^2+\frac{\gamma_c}{2}\int_{\Omega}\left(c(\mathbf{x},T)-\widehat{c}\right)^2+\frac{\gamma_u}{2}\int_{\partial\Omega\times(0,T)}u^2 \\
\ &\quad\quad+\frac{1}{2\epsilon}\int_{\Omega\times(0,T)}|\max\{0,u-u_+\}|^2+\frac{1}{2\epsilon}\int_{\Omega\times(0,T)}|\min\{0,u-u_-\}|^2,
\end{align*}
with $\epsilon$ a given (small) positive constant, chosen to enforce the control constraints efficiently.

When forming the Newton system in this setting, we will be required to solve systems relating to the finite element discretization of the following terms:
\begin{eqnarray}
\ \label{GN_ControlConstraints} \left[\begin{array}{ccccc}
\chi_{\Omega_T}(\sq) & 0 & 0 & \mathcal{L}_{zp} & \mathcal{L}_{zq} \\
0 & \gamma_{c}\chi_{\Omega_T}(\sq) & 0 & \mathcal{L}_{cp} & \mathcal{L}_{cq} \\
0 & 0 & \gamma_{u}\cdot\text{Id}+\frac{1}{\epsilon}G_{\Lambda} & 0 & -\beta\chi_{\partial\Omega}(\sq)^\top \\
\mathcal{L}_{pz} & \mathcal{L}_{pc} & 0 & 0 & 0 \\
\mathcal{L}_{qz} & \mathcal{L}_{qc} & -\beta\chi_{\partial\Omega}(\sq) & 0 & 0 \\
\end{array}\right]\left[\begin{array}{c}
s_z \\
s_c \\
s_u \\
s_p \\
s_q \\
\end{array}\right]=\widetilde{\mathbf{b}},
\end{eqnarray}
where
\begin{eqnarray*}
\ \widetilde{\mathbf{b}}:=\left[\begin{array}{c}
\widehat{z}-\left(\chi_{\Omega_T}(z)-\frac{\partial{}p}{\partial{}t}-D_{z}\nabla^{2}p-\alpha\nabla\left(\frac{1}{1+c}\right)\cdot\nabla{}p-2w\frac{zq}{(1+z^2)^2}\right) \\
\gamma_{c}\widehat{c}-\left(\gamma_{c}\chi_{\Omega_T}(c)-\alpha{}p\nabla\cdot\left(\nabla\left(\frac{1}{(1+c)^2}\right)z\right)-\frac{\partial{}q}{\partial{}t}-\nabla^{2}q+\rho{}q\right) \\
\frac{1}{\epsilon}(G_{\Lambda_+}y_{+}+G_{\Lambda_-}y_{-})-\left(\gamma_u{}u-\beta{}q\right) \\
-\left(\frac{\partial{}z}{\partial{}t}-D_{z}\nabla^{2}z-\alpha\nabla\cdot\left(\frac{\nabla{}c}{(1+c)^2}z\right)\right) \\
-\left(\frac{\partial{}c}{\partial{}t}-\nabla^{2}c+\rho{}c-w\frac{z^2}{1+z^2}\right) \\
\end{array}\right].
\end{eqnarray*}
Here, $G_{\Lambda_+}$, $G_{\Lambda_-}$, $G_{\Lambda}$ denote projections onto the active sets $\Lambda_+:=\{i:u_i>(u_+)_i\}$, $\Lambda_-:=\{i:u_i<(u_-)_i\}$, $\Lambda:=\Lambda_{+}\cup\Lambda_{-}$ (for the $i$-th node on the discrete level).

\section{Preconditioning for Gauss--Newton matrix systems}\label{sec:Preconditioner}

In this section we focus on deriving effective preconditioners for the matrix systems \eqref{GN} and \eqref{GN_ControlConstraints} resulting from the Gauss--Newton method applied to the chemotaxis problem, both without and with additional control constraints.

We base our preconditioners on the well studied field of \emph{saddle point systems}, which take the form \cite{BGL}
\begin{eqnarray}
\ \label{SaddlePt} \underbrace{\left[\begin{array}{cc}
A & B^\top \\
B & 0 \\
\end{array}\right]}_{\mathcal{A}}\left[\begin{array}{c}
\mathbf{x}_1 \\
\mathbf{x}_2 \\
\end{array}\right]=\left[\begin{array}{c}
\mathbf{b}_1 \\
\mathbf{b}_2 \\
\end{array}\right],
\end{eqnarray}
with $A$ symmetric positive semidefinite in our case, and $B$ having at least as many columns as rows. Two well-studied preconditioners for the system \eqref{SaddlePt} are given by \cite{Ipsen,Kuznetsov,MGW}
\begin{eqnarray*}
\ \mathcal{P}_D=\left[\begin{array}{cc}
A & 0 \\
0 & S \\
\end{array}\right],\quad\quad\mathcal{P}_T=\left[\begin{array}{cc}
A & 0 \\
B & -S \\
\end{array}\right],
\end{eqnarray*}
where the (negative) \emph{Schur complement} $S:=BA^{-1}B^\top$. It is known \cite{Ipsen,Kuznetsov,MGW} that, provided the preconditioned system is nonsingular, its eigenvalues are given by
\begin{equation*}
\ \lambda(\mathcal{P}_D^{-1}\mathcal{A})\in\left\{1,\frac{1}{2}(1\pm\sqrt{5})\right\},\quad\quad\lambda(\mathcal{P}_T^{-1}\mathcal{A})\in\left\{1\right\},
\end{equation*}
with these results also holding for the block triangular preconditioner $\mathcal{P}_T$ even if $A$ is not symmetric. Now, as $\mathcal{P}_D^{-1}\mathcal{A}$ is diagonalizable but $\mathcal{P}_T^{-1}\mathcal{A}$ is not, preconditioning with $\mathcal{P}_D$ ($\mathcal{P}_T$) yields convergence of a suitable Krylov subspace method in 3 (2) iterations, respectively.

In practice, however, $\mathcal{P}_D$ and $\mathcal{P}_T$ are not useful preconditioners, as the matrices $A$ and $S$ are computationally expensive to invert in general. We therefore instead seek preconditioners of the form
\begin{eqnarray*}
\ \widehat{\mathcal{P}}_D=\left[\begin{array}{cc}
\widehat{A} & 0 \\
0 & \widehat{S} \\
\end{array}\right],\quad\quad\widehat{\mathcal{P}}_T=\left[\begin{array}{cc}
\widehat{A} &0 \\
B & -\widehat{S} \\
\end{array}\right],
\end{eqnarray*}
where $\widehat{A}$ and $\widehat{S}$ denote suitably chosen approximations of the $(1,1)$-block $A$ and Schur complement $S$. The objective here is that our Krylov method will not converge in 3 or 2 iterations, but just a few more, while at the same time ensuring that our preconditioner is much cheaper to invert. From this point, we focus our attention on preconditioners of block triangular form.

\subsection{Construction of the preconditioner}\label{sec:Preconditioner_Construction}

We first examine the system \eqref{GN}, and place this in saddle point form \eqref{SaddlePt} as follows:
\begin{eqnarray*}
\ A=\left[\begin{array}{ccc}
\chi_{\Omega_T}(\sq) & 0 & 0 \\
0 & \gamma_{c}\chi_{\Omega_T}(\sq) & 0 \\
0 & 0 & \gamma_{u}\cdot\text{Id}
\end{array}\right],\quad\quad{}B=\left[\begin{array}{ccccc}
\mathcal{L}_{pz} & \mathcal{L}_{pc} & 0 \\
\mathcal{L}_{qz} & \mathcal{L}_{qc} & -\beta\chi_{\partial\Omega}(\sq) \\
\end{array}\right].
\end{eqnarray*}
Furthermore, let us decompose the blocks $A$ and $B$ into sub-blocks:
\begin{eqnarray*}
\ A=\left[\begin{array}{cc}
A_s & 0 \\
0 & A_u \\
\end{array}\right],\quad\quad{}B=\left[\begin{array}{cc}
B_s & B_u \\
\end{array}\right],
\end{eqnarray*}
where
\begin{eqnarray*}
\ A_s=\left[\begin{array}{cc}
\chi_{\Omega_T}(\sq) & 0 \\
0 & \gamma_{c}\chi_{\Omega_T}(\sq) \\
\end{array}\right],\quad\quad{}B_s=\left[\begin{array}{cc}
\mathcal{L}_{pz} & \mathcal{L}_{pc} \\
\mathcal{L}_{qz} & \mathcal{L}_{qc} \\
\end{array}\right],\quad\quad{}B_u=\left[\begin{array}{c}
0 \\
-\beta\chi_{\partial\Omega}(\sq) \\
\end{array}\right].
\end{eqnarray*}
In this paper, $A_u$ corresponds to the finite element discretization of the following operators:
\begin{eqnarray*}
\ A_u\leftarrow\left\{\begin{array}{cl}
\gamma_u\cdot\text{Id} & \text{without control constraints}, \\
\gamma_u\cdot\text{Id}+\frac{1}{\epsilon}G_{\Lambda} & \text{with control constraints}, \\
\end{array}\right.
\end{eqnarray*}
that is to say the block $A_u$ is altered if we instead consider the matrix \eqref{GN_ControlConstraints} incorporating control constraints. Note that, as the saddle point system is written, the matrix $A$ is not invertible and the Schur complement $S$ therefore does not exist. We hence consider a suitable re-ordering of the matrices \eqref{GN} and \eqref{GN_ControlConstraints} to enable us to utilize classical saddle point theory.

In particular, we observe that the matrix under consideration may be factorized as follows:
\begin{eqnarray*}
\ \left[\begin{array}{ccc}
A_s & 0 & B_s^\top \\
0 & A_u & B_u^\top \\
B_s & B_u & 0 \\
\end{array}\right]=\left[\begin{array}{ccc}
I & -A_{s}B_{s}^{-1}B_{u}A_u^{-1} & A_{s}B_s^{-1} \\
0 & I & 0 \\
0 & 0 & I \\
\end{array}\right]\underbrace{\left[\begin{array}{ccc}
0 & 0 & S_{\angle} \\
0 & A_u & B_u^\top \\
B_s & B_u & 0 \\
\end{array}\right]}_{\mathcal{P}},
\end{eqnarray*}
where identity matrices $I$ are of appropriate dimensions. Note that as $\text{Id}$ corresponds to the identity operator on the continuous level, this will become a finite element mass matrix in the discrete setting. We then take $\mathcal{P}$ to be the foundation of our preconditioner. We define $S_{\angle}$ as the `\emph{pivoted Schur complement}' \cite[Section 3.3]{DPSS}
\begin{eqnarray}
\ \label{Schur} S_{\angle}=B_s^\top+A_{s}B_{s}^{-1}B_{u}A_{u}^{-1}B_u^\top.
\end{eqnarray}
We approximate this term within our preconditioner using the `\emph{matching strategy}' devised in \cite{PSW,PW2,PW1}, which aims to capture both terms of the Schur complement within the preconditioner. The approximation reads as follows:
\begin{eqnarray}
\ \label{Sangle} S_{\angle}\approx\widehat{S}_{\angle}:=\Big(B_s^\top+\frac{1}{\eta}A_s\Big)B_s^{-1}\Big(B_s+\eta{}B_{u}A_u^{-1}B_u^\top\Big),
\end{eqnarray}
Note that the matrix product $B_s^\top{}B_s^{-1}B_s$ captures the first term $B_s^\top$ of $S_{\angle}$, and $\big(\frac{1}{\eta}A_s\big)B_s^{-1}\big(\eta{}B_{u}A_u^{-1}B_u^\top\big)$ matches exactly the second term $A_{s}B_{s}^{-1}B_{u}A_{u}^{-1}B_u^\top$. The positive constant $\eta$ is chosen to `balance' the first and last matrix factors, $B_s^\top+\frac{1}{\eta}A_s$ and $B_s+\eta{}B_{u}A_u^{-1}B_u^\top$, within the Schur complement approximation, so that the two terms in the remainder $S_{\angle}-\widehat{S}_{\angle}$ are approximately of the same norm.
Two natural choices for this constant are
\begin{eqnarray*}
\ \eta=\sqrt{\frac{\left\|A_s\right\|}{\left\|B_{u}A_u^{-1}B_u^\top\right\|}}\quad\text{or}\quad\eta=\sqrt{\frac{\max(\text{diag}(A_s))}{\max(\text{diag}(B_{u}A_u^{-1}B_u^\top))}},
\end{eqnarray*}
with the second such choice much cheaper to compute. Approximately solving for the matrix $B_s+\eta{}B_{u}A_u^{-1}B_u^\top$ is made tractable by the effective approximation of a mass matrix (or a mass matrix plus a positive diagonal matrix) of the form $A_u$ by its diagonal, see \cite[Section 4.1]{PearsonGondzio} and \cite{WathenEigBounds}.

Putting all the pieces together, we state our preconditioner
\begin{eqnarray*}
\ \widehat{\mathcal{P}}=\left[\begin{array}{ccc}
0 & 0 & \widehat{S}_{\angle} \\
0 & A_u & B_u^\top \\
B_s & B_u & 0 \\
\end{array}\right],
\end{eqnarray*}
incorporating the Schur complement approximation above. Due to the re-ordering of the saddle point system that we have undertaken, this is a suitable choice of preconditioner that captures the characteristics of the matrix under consideration.

\subsection{Application of the preconditioner}\label{sec:Preconditioner_Application}

Applying the inverse of the preconditioner, $\widehat{\mathcal{P}}^{-1}$, as is necessary within an iterative method, therefore requires three main operations:
\begin{enumerate}
	\item \underline{Applying $B_s^{-1}$:} This is equivalent to solving the forward problem, rather than the coupled optimization problem. In practice this is approached time-step by time-step, using an algebraic or geometric multigrid method, or another suitable scheme, to solve for the matrices arising at each point in time.

	\item \underline{Applying $A_u^{-1}$:} The matrix $A_u$ is a block diagonal matrix, consisting of boundary mass matrices at each time-step (in the case without control constraints), or boundary mass matrices plus positive semidefinite diagonal matrices (if control constraints are present). In either case, these matrices may be well approximated using Chebyshev semi-iteration \cite{GVI,GVII,WathenRees}, or even using a simple diagonal approximation of a mass matrix \cite{WathenEigBounds}.

	\item \underline{Applying $\widehat{S}_{\angle}^{-1}$:} Applying the approximation \eqref{Sangle} involves a multiplication operation involving $B_s$, and (approximate) solves for each of $B_s^\top+\frac{1}{\eta}A_s$ and $B_s+\eta{}B_{u}A_u^{-1}B_u^\top$ which may again be approached at each time-step in turn using multigrid or another appropriate method.
\end{enumerate}

\subsection{Uzawa approximation}\label{sec:Preconditioner_Uzawa}

In practice, we make a further modification to the preconditioner $\widehat{\mathcal{P}}$ in order to ensure it is easier to work with on a computer. In more detail, the term $B_s$ in the bottom-left of $\widehat{\mathcal{P}}$, and the terms $B_s^\top+\frac{1}{\eta}A_s$ and $B_s+\eta{}B_{u}A_u^{-1}B_u^\top$ within $\widehat{S}_{\angle}$, contain $2\times2$ block systems which we would like to replace with more convenient approximations so that we are only required to (approximately) invert one block at a time.

To facilitate this, we replace $2\times2$ block matrices by an inexact Uzawa approximation, with block triangular splitting matrices, where appropriate. This leads to our final choice of preconditioner:
\begin{eqnarray*}
\ \widehat{\mathcal{P}}_{\text{Uzawa}}=\left[\begin{array}{ccc}
0 & 0 & \Big(B_s^\top+\frac{1}{\eta}A_s\Big)_{\text{Uzawa}}(B_s)_{\text{Uzawa}}^{-1}\Big(B_s+\eta{}B_{u}A_u^{-1}B_u^\top\Big)_{\text{Uzawa}} \\
0 & A_u & B_u^\top \\
(B_s)_{\text{Uzawa}} & B_u & 0 \\
\end{array}\right],
\end{eqnarray*}
where $(\cdot)_{\text{Uzawa}}$ denotes the Uzawa approximation of the corresponding matrix.
For ease of reproducibility for the reader, we state the splitting matrices below:
\begin{align*}
\ (B_s)_{\text{Uzawa}}\rightarrow{}&\left[\begin{array}{cc}
\mathcal{L}_{pz} & \mathcal{L}_{pc} \\
0 & \mathcal{L}_{qc} \\
\end{array}\right], \\
\ \Big(B_s^{\top}+\frac{1}{\eta}A_s\Big)_{\text{Uzawa}}\rightarrow{}&\left[\begin{array}{cc}
\mathcal{L}_{zp}+\frac{1}{\eta}\chi_{\Omega_T}(\sq) & 0 \\
\mathcal{L}_{cp} & \mathcal{L}_{cq}+\frac{\gamma_c}{\eta}\chi_{\Omega_T}(\sq) \\
\end{array}\right], \\
\ \Big(B_s+\eta{}B_{u}A_u^{-1}B_u^\top\Big)_{\text{Uzawa}}\rightarrow{}&\left[\begin{array}{cc}
\mathcal{L}_{pz} & \mathcal{L}_{pc} \\
0 & \mathcal{L}_{qc}+\eta\beta^2\chi_{\partial\Omega}(\sq)A_u^{-1}\chi_{\partial\Omega}(\sq)^\top \\
\end{array}\right].
\end{align*}
Now the linear systems with diagonal blocks ($\mathcal{L}_{zp}$, $\mathcal{L}_{qc}$, and so on) can be solved directly.
Note that it is also possible to annihilate another off-diagonal block instead within the Uzawa approximation.
However, we have found that the approximations listed above yield fast convergence in the numerical experiments.

\section{Numerical experiments with control constraints}\label{sec:NumEx1}
In this section we benchmark the preconditioned Newton method.
For our test problem, the initial distribution of bacterial cells is chosen as a sum of $m_0$ independent Gaussian peaks,
\begin{equation}
 z_0(x,y) = \sum_{i=1}^{m_0} \exp\left(-2560 \cdot \left[(x-x_i)^2 +(y-y_i)^2 \right]\right),
\label{eq:z0_full}
\end{equation}
where the centers $\{x_i,y_i\}$ are chosen randomly on $[0,1]^2$.
The desired distribution at the final time $T=1$ is linear,
\begin{equation}
\widehat z(x,y) = \langle z_0 \rangle \cdot (x+y),
\label{eq:zhat}
\end{equation}
normalized by the initial mass,
$$
\langle z_0 \rangle = \int_{[0,1]^2} z_0(x,y)~{\rm d}x{\rm d}y,
$$
since the model conserves the normalization of $z$.
Both initial and target concentrations $c$ are zero.
The experiments were run in {\scshape matlab} R2017b on one core of a 2.4GHz Intel Xeon E5-2640 CPU.

In this section, we set $m_0=50$ and the control constraints $u_-=0$ and $u_+=0.2$, in accordance with \cite{Potschka}.
The default regularization parameters are set to $\gamma_u=10^{-3}$ and $\gamma_c=0.5$.
The stopping tolerance for the Newton iteration is set to $10^{-4}$.
Moreover, we decrease the Moreau--Yosida regularization parameter $\epsilon$ geometrically from $10^{-1}$ to $10^{-4}$ as the iteration converges.
This gives more robust behavior of the Newton method.

%
%
%
%

\begin{figure}[htb]
\centering
\caption{Left: CPU time of the Newton solver for the constrained control. Right: $u({\bf x},T/2)$.}
\label{fig:full_ttimes}
\resizebox{0.49\linewidth}{!}{%
\begin{tikzpicture}
\begin{axis}[%
ylabel=CPU time, axis y line=left, y label style={at={(-0.05,1.0)},anchor=south west,rotate=-90}, ymode=log, yminorgrids,%
xlabel=$\log_{2}n$,%
legend style={at={(0.01,0.99)},anchor=north west}]

\addplot+[blue,solid] coordinates{(5,146.6925)
                                  (6,1508.9)
                                  (7,14186)
                                 }; \addlegendentry{Measured};
\addplot+[blue,dashed,domain=5:8,no marks] {2^(3.3*x-9.3)}; \addlegendentry{$n^{3.3}$};
\end{axis}
\end{tikzpicture}
}
\includegraphics[width=0.49\linewidth]{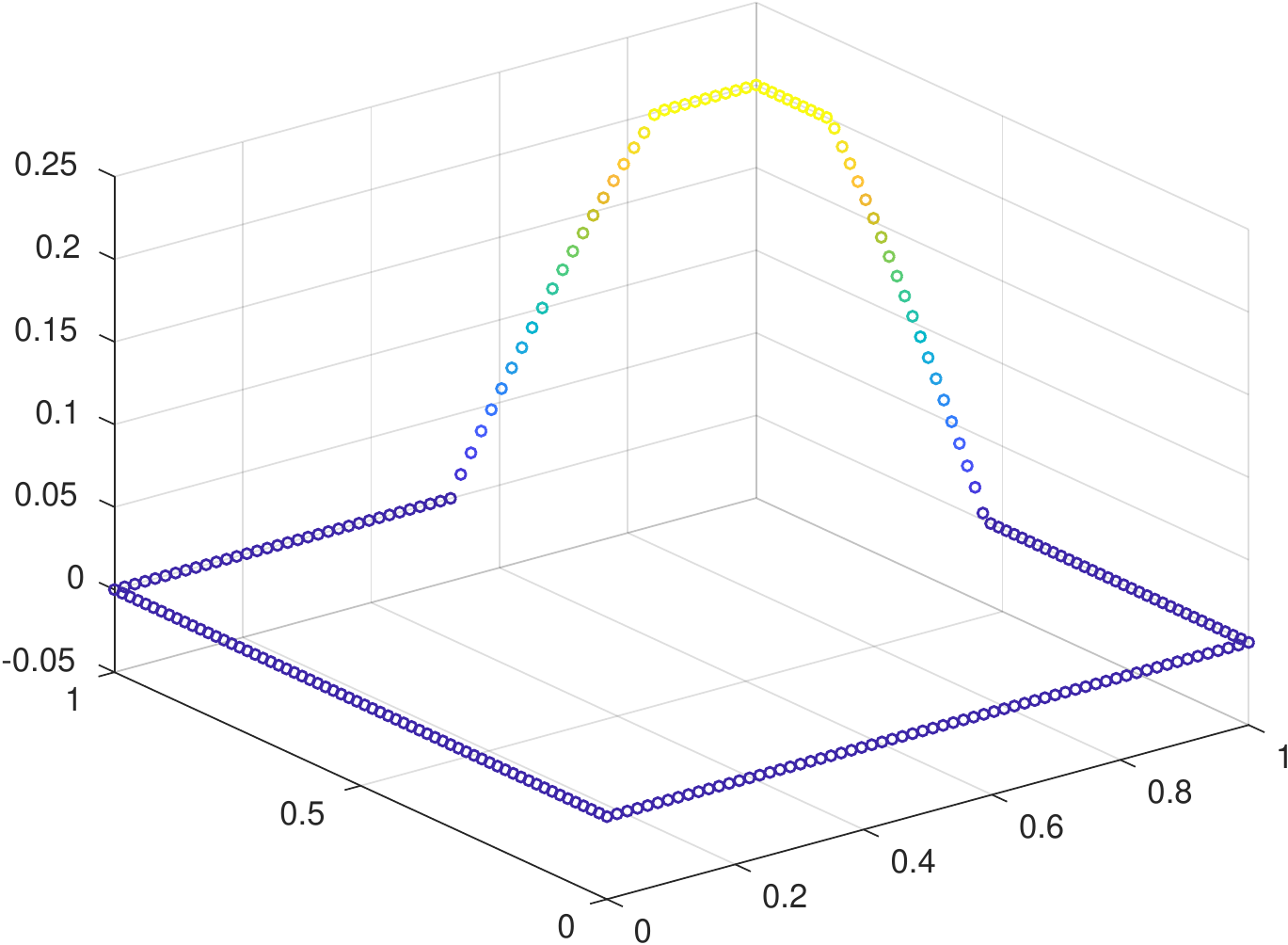}
%
\end{figure}
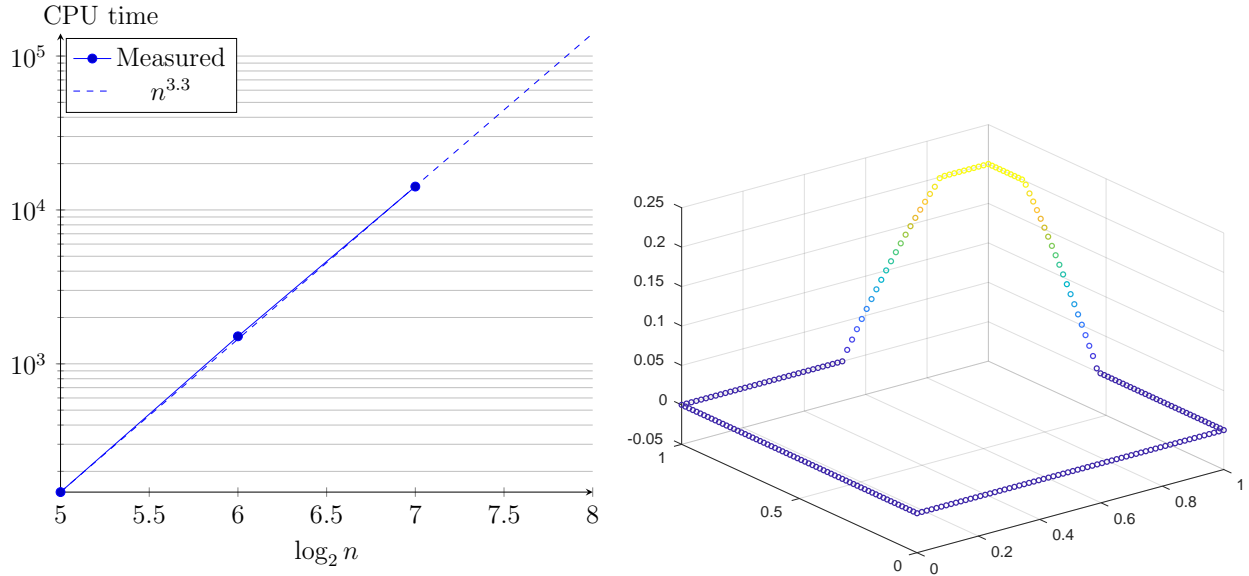

\begin{figure}[htb]
\centering
\caption{Left: cell density at the final time $z({\bf x},T)$. Right: misfit of the density, $z({\bf x},T) - \widehat z({\bf x})$.}
\label{fig:full_Z}
\includegraphics[width=0.49\linewidth]{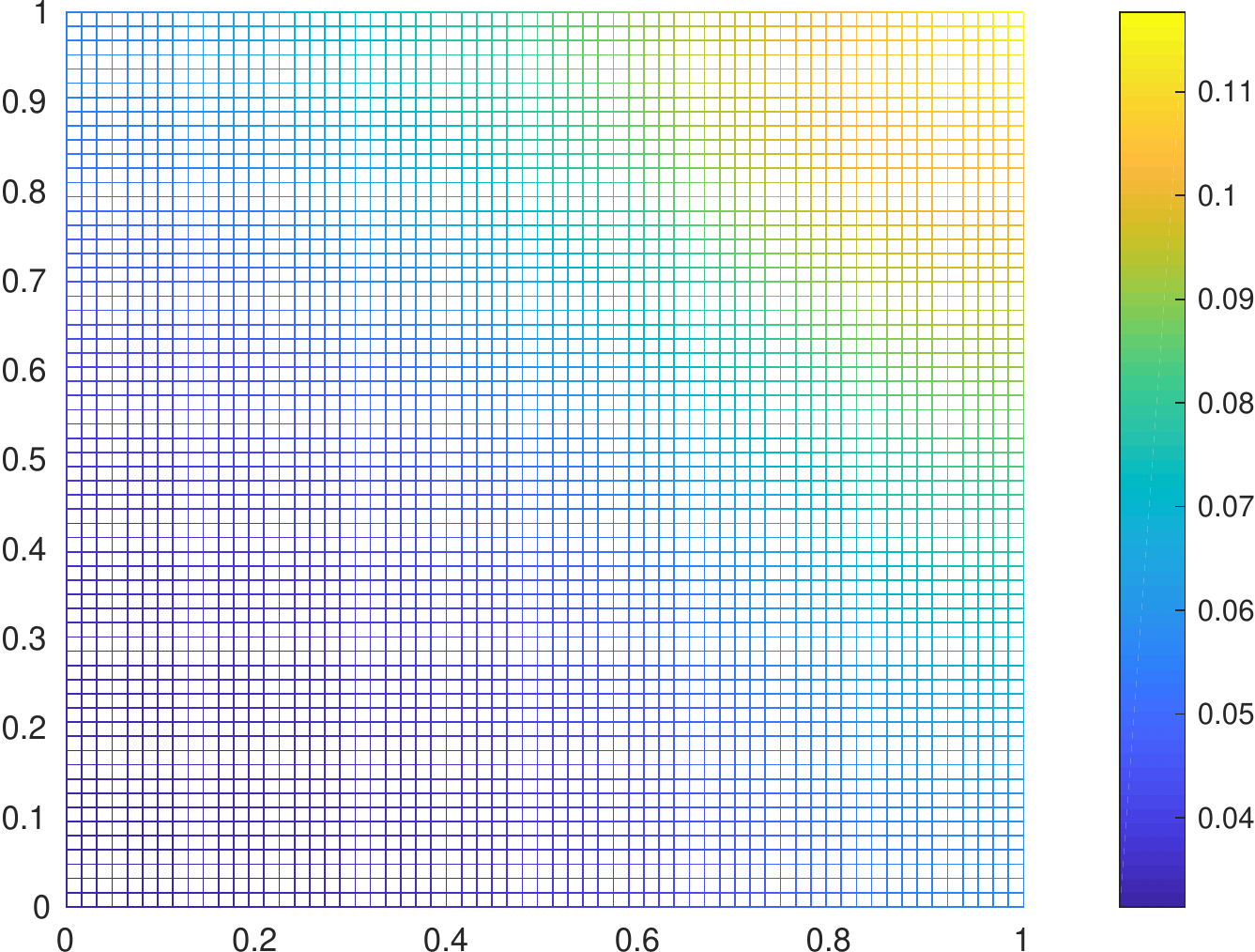}
\includegraphics[width=0.49\linewidth]{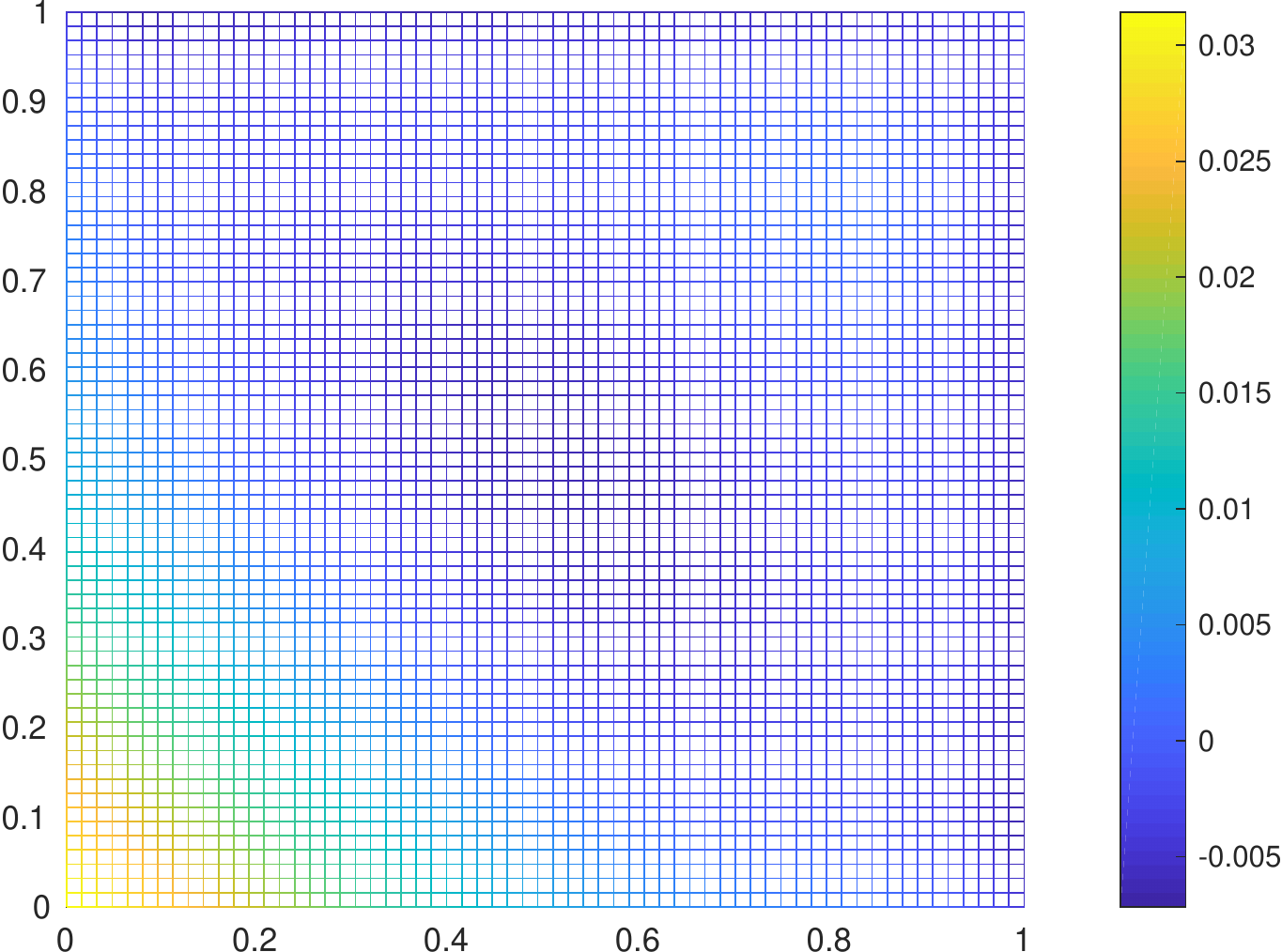}
\end{figure}

The computational time is shown in Fig. \ref{fig:full_ttimes} (left).
We see that it grows cubically with respect to the uniform grid refinement, which is expected for a three-dimensional (2D space + time) problem.
The number of Newton iterations is quite stable with respect to the grid size, ranging from $11$ to $14$ depending on a particular distribution of the random initial guess.

The transient control signal is shown in Fig. \ref{fig:full_ttimes} (right).
We notice that it is accurately confined within the prescribed constraints.
However, this leads to a rather large misfit in the target cell density (Fig. \ref{fig:full_Z}).
While the density follows the linear distribution $\widehat z$ correctly in the top right corner of the domain,
in the left bottom corner we see an excessive density of bacteria.
This shows that controlling only the chemoattractant might be insufficient for forcing the bacteria to leave a particular area.

Lastly in this section, we investigate the performance of the preconditioner proposed in Section \ref{sec:Preconditioner} against variation of parameters.
In Table \ref{tab:full_its}, we show the average numbers of GMRES \cite{gmres} iterations per Newton step for different grid sizes $n$ and regularization parameters $\gamma_u$, $\gamma_c$.
We vary only one parameter at a time, while the other two are kept fixed to their default values, $n=64$, $\gamma_u=10^{-3}$, and $\gamma_c=0.5$.
The number of iterations grows slightly as the control regularization parameter $\gamma_u$ is decreased, which is expected for a boundary control problem.
On the other hand, the preconditioner is reasonably robust with respect to the other parameters, in particular the grid size.

\begin{table}[htb]
 \centering
 \caption{Average number of GMRES iterations per Newton step.}
 \label{tab:full_its}
 \begin{tabular}[t]{c|c}
 $n$  & its \\ \hline
 32   & 21.37  \\
 64   & 27.46  \\

 128  & 27.86  \\
 \end{tabular}\hfill
 \begin{tabular}[t]{c|c}
 $\gamma_u$  &  its  \\ \hline
 $10^0$      &  6.00        \\
 $10^{-1}$   &  9.00        \\
 $10^{-2}$   &  15.28  \\
 $10^{-3}$   &  27.46  \\
 $10^{-4}$   &  46.84  \\
 $10^{-5}$   &  69.21  \\
 \end{tabular}\hfill
 \begin{tabular}[t]{c|c}
 $\gamma_c$  &  its  \\ \hline
 $0.5$       &  27.46   \\
 $10^{-1}$   &  30.55   \\
 $10^{-2}$   &  32.11   \\
 $10^{-3}$   &  31.77   \\
 $10^{-4}$   &  32.67   \\
 $10^{-5}$   &  31.57   \\
 \end{tabular}
\end{table}

\section{Low-rank tensor decompositions and algorithms} \label{sec:LowRank}
The optimality system \eqref{GN} can result in a huge-scale matrix system, for many spatial degrees of freedom and time steps.
One way to reduce the associated computational burden is to seek an approximate solution in a low-parametric representation.
In this paper we apply separation of variables, and in particular the Tensor Train (TT) decomposition \cite{osel-tt-2011}.
In this section, we introduce the TT decomposition and the algorithm for an efficient TT-structured solution of the optimality equations.
Although the TT approximation can have difficulties with the indicator function of the active set of control constraints (see Remark \ref{rem:indicator} below),
for the problem without box constraints it yields a very efficient solver.
So in this section we assume an unconstrained control setting.

\subsection{Tensor product discretization and indexing}
We assume that the solution functions can be discretized on a structured grid, e.g. the cell concentration $z({\bf x},t)$ with a $d$-dimensional spatial variable ${\bf x} = (x_1,\ldots,x_d)$ can be approximated by
$$
z({\bf x},t) \approx \sum_{i_1,\ldots,i_d, i_{d+1}=1}^{n_1,\ldots,n_d, n_{d+1}} {\bf z}(i_1,\ldots,i_d,i_{d+1}) \phi_{i_1,\ldots,i_d}({\bf x}) \psi_{i_{d+1}}(t),
$$
where $\{\phi_{i_1,\ldots,i_d}({\bf x})\}$ is a set of spatial basis functions as introduced in Section \ref{sec:Matrix}, which we now assume to be indexed by $d$ independent variables.
In particular, we consider a square domain ${\bf x} \in [0,1]^d$ and the piecewise polylinear basis functions
$$
\phi_{i_1,\ldots,i_d}({\bf x}) = \varphi_{i_1}(x_1) \cdots \varphi_{i_d}(x_d).
$$
In turn, $\{\psi_{i_{d+1}}(t)\}$ is a set of nodal interpolation functions in time, associated with the uniform time grid $\{t_{i_{d+1}}\}$, with $t_{i_{d+1}} = \tau \cdot i_{d+1}$, $i_{d+1} = 1,\ldots,n_{d+1}$, and $\tau=T/n_{d+1}$.
We can see that the discrete coefficients of $z$ can be collected into a $(d+1)$-dimensional \emph{tensor}.
Introducing a uniform bound $n \ge n_k$, $k=1,\ldots,d+1$, we can immediately conclude that the tensor ${\bf z}$ has $\mathcal{O}(n^{d+1})$ entries.
The computational complexity of solving \eqref{GN} is usually much higher.
This explains the sometimes relatively high computing times in the previous section.

\emph{Separation} of the discrete variables $i_1,\ldots,i_{d+1}$ can compress the tensor data from the exponential $\mathcal{O}(n^{d+1})$ to a linear volume $\mathcal{O}(dn)$.
Yet we can aim for a higher compression ratio.
Assuming that the range $n_k$ of an index $i_k$ is factorizable into a set of divisors $n_{k,1}\cdots n_{k,L_k}=n_k$,
we can also factorize the index $i_k$ into the corresponding digits,
$$
i_k = 1 + \sum_{\ell=1}^{L_k} (i_{k,\ell}-1) \prod_{p=1}^{\ell-1} n_{k,p}, \qquad k=1,\ldots,d+1.
$$
Now the tensor ${\bf z}$ can be enumerated by the elementary digits $i_{k,\ell}$, which we shall denote simply as $i_m$ from now on, for $m=1,\ldots,L = \sum_{k=1}^{d+1} L_k$.
Instead of considering ${\bf z}$ as a $(d+1)$-dimensional tensor, we treat it as a $L$-dimensional tensor with elements ${\bf z}(i_1,\ldots,i_L)$, and therefore we will separate now the \emph{virtual} indices $i_m$ \cite{tee-tensor-2003}.

\subsection{Tensor Train decomposition}
As a particular separated approximation, we choose the \emph{Tensor Train} (TT) decomposition \cite{osel-tt-2011}, which is also known as the Matrix Product States \cite{PerezGarcia-mps-2007,schollwock-2005} in physics:
\begin{equation}
{\bf z}(i_1,\ldots,i_L) \approx \sum\limits_{s_1=1}^{r_1} \cdots \sum\limits_{s_{L-1}=1}^{r_{L-1}} z^{(1)}_{s_1}(i_1) z^{(2)}_{s_1,s_2}(i_2) \cdots z^{(L)}_{s_{L-1}}(i_L).
\label{eq:tt}
\end{equation}
The factors $z^{(m)}$ on the right hand side are called \emph{TT blocks}, and the ranges $r_1,\ldots,r_{L-1}$ of the auxiliary summation indices are called \emph{TT ranks}.
Notice that the TT blocks are at most $3$-dimensional tensors, of sizes $r_{m-1} \times n_m \times r_m$ (for uniformity, we can let $r_0=r_L=1$).

Potentially, we can represent any finite dimensional tensor exactly through \eqref{eq:tt} by choosing large enough TT ranks.
For reasons of numerical efficiency, we will of course aim for a (sub-)optimal approximation with $r_m$ being as small as possible, and most importantly much smaller than the original tensor size $n_1\cdots n_{d+1}$.
The storage needed for the right hand side of \eqref{eq:tt} is of the order of $\mathcal{O}(L n_m r_m^2)$,
where $n_m=n_{k,\ell}$ is also chosen to be much smaller than the original $n_k$.
For example, if we restrict the grid sizes to be powers of two, $n_k=2^{L_k}$, the range of each index $i_m$ in \eqref{eq:tt} becomes just $\{1,2\}$, whereas $L$, and hence the storage complexity of the TT format, becomes \emph{logarithmic} in the original tensor size, $L = \log_2(n_1\cdots n_{d+1})$.
Due to the minimal non-trivial index range in this case, the TT decomposition \eqref{eq:tt} with $i_m \in \{1,2\}$ was called the \emph{Quantized} TT (QTT) decomposition \cite{khor-qtt-2011}.
It was then proved that many examples of vectors \cite{khor-qtt-2011} and matrices \cite{khkaz-lap-2012,khkaz-conv-2013}, arising from the discretization of functions and differential operators, allow low-rank QTT decompositions.

Abstracting from the original problem dimensions, we can consider only two data representations: a tensor with the smallest possible ranges ${\bf z}(i_1,\ldots,i_L)$, and a vector of the same data entries:
\begin{equation}
{\bf z}(i) = {\bf z}(i_1,\ldots,i_L), \quad \mbox{where} \quad i = 1 + \sum_{m=1}^{L} (i_{m}-1) \prod_{p=1}^{m-1} n_p.
\label{eq:ten-vec}
\end{equation}
We need the vector notation for setting the Gauss--Newton equations \eqref{GN} on tensors consistently.
Boldface letters (e.g. ${\bf z}$) from now on will denote vectors.
We can use the Kronecker product ($\otimes$) to rewrite \eqref{eq:tt} in an equivalent vector form,
\begin{equation*}
{\bf z} = \sum\limits_{s_1,\ldots,s_{L-1}=1}^{r_1,\ldots,r_{L-1}} z^{(1)}_{s_1} \otimes z^{(2)}_{s_1,s_2} \otimes \cdots \otimes  z^{(L)}_{s_{L-1}}.
\end{equation*}
Of course, we shall never actually compute the Kronecker products in the expansion above, but only store and manipulate individual TT blocks on the right hand side.

For example, the matrix-vector product ${\bf y} = A{\bf z}$ with ${\bf z}$ given in \eqref{eq:tt} can be computed efficiently if we can also represent the matrix by a TT decomposition,
\begin{equation}
A = \sum\limits_{s_1,\ldots,s_{L-1}=1}^{R_1,\ldots,R_{L-1}} A^{(1)}_{s_1} \otimes A^{(2)}_{s_1,s_2} \otimes \cdots \otimes  A^{(L)}_{s_{L-1}}.
\label{eq:ttm}
\end{equation}
For example if the matrix is diagonal, and the vector of the diagonal values can be represented by a TT decomposition \eqref{eq:tt}, the matrix can be written as in \eqref{eq:ttm}, with the same TT ranks.
There are less trivial matrices, arising for example in finite element computations, that admit TT decompositions with modest ranks $R_m$ \cite{khkaz-lap-2012,khkaz-conv-2013}.
Now the result ${\bf y}=A{\bf z}$ can be also written in the TT format and computed block by block.
Moreover, a TT decomposition with excessive TT ranks can be efficiently approximated up to a desired accuracy by a decomposition with sub-optimal ranks using QR and singular value decomposition (SVD) factorizations \cite{osel-tt-2011}, without ever constructing full large tensors.

\subsection{Alternating Linear Scheme iteration for solving \eqref{GN}}

In addition to the cell concentration $z({\bf x},t)$, we need to represent the other solution components.
Since all components are defined on the same domain, we can discretize them using the same basis.
The tensors of discrete values therefore have the same sizes.
The structure of the problem \eqref{GN} suggests that we approximate them in a shared TT decomposition, the so-called \emph{block} TT format \cite{dkos-eigb-2014}.
We denote the aggregated solution
\begin{equation*}
{\bf y}^\top = \begin{bmatrix}
           {\bf z}^\top & {\bf c}^\top & {\bf p}^\top & {\bf q}^\top & {\bf u}^\top
          \end{bmatrix},
\end{equation*}
enumerating the components via ${\bf y}_{j}$, $j=1,\ldots,5$.
Now we decompose ${\bf y}$ into a TT format with all the same TT blocks except the $m$-th block for some $m=1,\ldots,L$, which actually carries the enumerator of the components,
\begin{equation}
{\bf y}_{j} = \sum_{s_1,\ldots,s_{L-1}=1}^{r_1,\ldots,r_{L-1}} y^{(1)}_{s_1} \otimes \cdots \otimes y^{(m-1)}_{s_{\ell-2},s_{\ell-1}} \otimes \widehat y^{(m)}_{s_{\ell-1},s_{\ell}}(j) \otimes y^{(m+1)}_{s_{\ell},s_{\ell+1}} \otimes  \cdots \otimes y^{(L)}_{s_{L-1}}.
\label{eq:btt}
\end{equation}
Moreover, we can switch between the representations \eqref{eq:btt} corresponding to different $m$ (and hence having $j$ in different TT blocks) using the SVD \cite{dkos-eigb-2014}.
For example, we can reshape $\widehat y^{(m)}$ into a matrix with elements
$$
\widehat Y^{(m)}(s_{m-1},i_m;~j,s_m) = \widehat y^{(m)}_{s_{m-1},s_m}(i_m,j)
$$
and compute the truncated SVD $\widehat Y^{(m)} \approx U \Sigma V^\top$.
Now we write the left singular vectors $U$ into the $m$-th TT block instead of $\widehat y^{(m)}$, and multiply $\Sigma V^\top$ with the $(m+1)$-th TT block,
\begin{align}
\label{eq:svd1} y^{(m)}_{s_{m-1},s_m'}(i_m)={}&U(s_{m-1},i_m;~s_m'), \\
\label{eq:svd2} \widehat y^{(m+1)}_{s_m',s_{m+1}}(i_{m+1},j)={}&\sum_{s_m=1}^{r_m}\Sigma V^\top (s_m';~j,s_m) y^{(m+1)}_{s_m,s_{m+1}}(i_{m+1}).
\end{align}
Note that we have obtained the same representation as \eqref{eq:btt} with $m$ replaced by $m+1$.
This process can be continued further, or reversed, and hence the $j$-index can be placed into any TT block.

A crucial ingredient for the iterative computation of \eqref{eq:btt} is the \emph{linearity} of the TT format.
Having chosen an $m=1,\ldots,L$, we construct the so-called \emph{frame} matrix, where the TT block $\widehat y^{(m)}$ in \eqref{eq:btt} is replaced by the identity matrix,
\begin{equation}
  Y_{m} = \sum_{s_1,\ldots,s_{m-2}}y^{(1)}_{s_1}\otimes \cdots \otimes y^{(m-1)}_{s_{m-2}} \otimes I_{n_m} \otimes \sum_{s_{m+1},\ldots,s_{L-1}} y^{(m+1)}_{s_{m+1}} \otimes \cdots \otimes y^{(L)}_{s_{L-1}}.
 \label{eq:frame}
\end{equation}
If we now treat $\widehat y^{(m)}(j)$ as a vector, we can observe that
\begin{equation}
{\bf y}_j = Y_{m} \widehat y^{(m)}(j),
\label{eq:ttlin}
\end{equation}
i.e. the frame matrix realises a linear map from the elements of $\widehat y^{(m)}$ to the elements of the whole solution vectors.
This motivates an iterative algorithm \cite{holtz-ALS-DMRG-2012}, which was called the Alternating Linear Scheme (ALS):

\begin{algorithmic}[1]
 \For {$\mbox{iter}=0,1,\ldots$ until convergence}
   \For {$m=1,2,\ldots,d,d-1,\ldots,1$}
     \State Plug the solution in the form \eqref{eq:ttlin} into the original problem.
     \State Solve the resulting overdetermined problem on $\widehat y^{(m)}$.
     \State Prepare the format \eqref{eq:btt} and the frame matrix \eqref{eq:frame} for $m+1$ or $m-1$.
   \EndFor
 \EndFor
\end{algorithmic}

Starting from a low-rank initial guess of the form \eqref{eq:btt}, this algorithm seeks the solution in a low-rank TT format by sweeping through the different TT blocks.
However, there might be different ways to resolve the overdetermined problem in Line 4.
For the optimality equations of the inverse problem, such as in \eqref{GN}, it was found \cite{bdos-sb-2016,ds-navier-2017} to be efficient to
use columns of the frame matrix as a Galerkin basis and project each submatrix of the Karush--Kuhn--Tucker (KKT) system individually.
In our case we notice that the $(3,3)$-block of \eqref{GN}
is simply a diagonal matrix in the case where lumped mass matrices are considered, and therefore eliminate the control component from the equations.\footnote{Our derivation is of course valid for any invertible matrix $A_u$, however we wish to exploit the simplicity of the matrix structure within our solver. When consistent mass matrices are applied, we can well approximate these by their diagonals within a preconditioner, see \cite{WathenEigBounds}.}
Specifically, we deduce that $s_u = A_u^{-1} \left({\bf b}_u + \beta\chi_{\partial\Omega}^\top s_q\right)$ and plug this into the fifth row.
This gives us a system of 4 equations only.
Moreover, instead of using the increments $s_z,s_c,s_p,s_q$, we can rewrite the equations for the new solution components directly:
$$
\begin{bmatrix}
 \chi_{\Omega_T} & 0 & \mathcal{L}_{zp} & \mathcal{L}_{zq} \\
 0 & \gamma_{c}\chi_{\Omega_T} & \mathcal{L}_{cp} & \mathcal{L}_{cq} \\
 \mathcal{L}_{pz} & \mathcal{L}_{pc} & 0 & 0 \\
 \mathcal{L}_{qz} & \mathcal{L}_{qc} & 0 & -\beta^2\chi_{\partial\Omega}A_u^{-1}\chi_{\partial\Omega}^\top \\
\end{bmatrix}
\begin{bmatrix}
{\bf z} \\
{\bf c} \\
{\bf p} \\
{\bf q}
\end{bmatrix}
=
\begin{bmatrix}
{\bf \widetilde b}_z \\
{\bf \widetilde b}_c \\
{\bf \widetilde b}_p \\
{\bf \widetilde b}_q
\end{bmatrix},
$$
where ${\bf \widetilde b}$ is the correspondingly adjusted right hand side.
Now we plug in the solutions in the form \eqref{eq:ttlin} (with $j$ now running only from $1$ to $4$), and project each of the previous equations onto $Y_m$.
This gives us a \emph{reduced} system
\begin{equation}
 \begin{bmatrix}
 \widehat \chi_{\Omega_T} & 0 & \mathcal{\widehat L}_{zp} & \mathcal{\widehat L}_{zq} \\
 0 & \gamma_{c}\widehat \chi_{\Omega_T} & \mathcal{\widehat L}_{cp} & \mathcal{\widehat L}_{cq} \\
 \mathcal{\widehat L}_{pz} & \mathcal{\widehat L}_{pc} & 0 & 0 \\
 \mathcal{\widehat L}_{qz} & \mathcal{\widehat L}_{qc} & 0 & -\beta^2\widehat \chi^2_{\partial\Omega} \\
\end{bmatrix}
\widehat y^{(m)}
=
\begin{bmatrix}
Y_m^\top {\bf \widetilde b}_z \\
Y_m^\top {\bf \widetilde b}_c \\
Y_m^\top {\bf \widetilde b}_p \\
Y_m^\top {\bf \widetilde b}_q
\end{bmatrix},
\label{eq:GN-red}
\end{equation}
with $\widehat \chi_{\Omega_T} = Y_m^\top \chi_{\Omega_T} Y_m$, $\mathcal{\widehat L}_{**} = Y_m^\top \mathcal{L}_{**} Y_m$ (where ``$**$'' stands for ``$zp$'', ``$zq$'', ``$cp$'', and so on), and $\widehat \chi^2_{\partial\Omega} = Y_m^\top \chi_{\partial\Omega}A_u^{-1}\chi_{\partial\Omega}^\top Y_m$ the projected square matrices.
Each submatrix is of size $r_{m-1}n_m r_m \times r_{m-1}n_m r_m$ (remember that we can choose $n_m=2$), and hence \eqref{eq:GN-red} is easy to solve.
Moreover, the singular value decomposition in \eqref{eq:svd1}--\eqref{eq:svd2} maintains the orthogonality of the frame matrices $Y_m$ automatically in the course of alternating iterations, provided that the initial guess is given with this property.
This makes the projected submatrices well conditioned if the original matrices were so, which eventually makes the entire matrix in \eqref{eq:GN-red} invertible.
We highlight that the preconditioner developed in Section \ref{sec:Preconditioner} can also be used for solving the system \eqref{eq:GN-red}.

\subsection{Construction of matrices in the TT format}
In the course of the Newton iteration, we need to reconstruct the matrices in \eqref{GN} (and consequently in \eqref{eq:GN-red}) using the new solution.
Assume that we need to construct an abstract bilinear form of a nonlinear transformation $f$ of the solution,
\begin{equation}
\mathcal{L}_{\mathcal{B}} = \int f(z,c,p,q) \nabla^{p} \phi_i \cdot \nabla^q \phi_j~{\rm d}{\bf x},
\label{eq:bilin-abs}
\end{equation}
where $p,q \in \{0,1\}$ are the differentiation orders, and $\phi_i,\phi_j$ are the basis functions.
Instead of the exact functions $z,c,p,q$, we work with the tensors of their values, ${\bf z},{\bf c},{\bf p},{\bf q}$.
The corresponding values of $f$ can be also collected into a tensor ${\bf f}$ of the same size, and the original function can be approximated in the same basis, i.e.
\begin{equation*}
f(z({\bf x}),c({\bf x}),p({\bf x}),q({\bf x})) \approx \sum_{i_1,\ldots,i_d} {\bf f}(i_1,\ldots,i_d) \phi_{i_1,\ldots,i_d}({\bf x}).
\end{equation*}
Now the computation of \eqref{eq:bilin-abs} involves computing analytical triple products
\begin{equation}
 \mathcal{H}(i,j,k) = \int \phi_k \nabla^{p} \phi_i \cdot \nabla^q \phi_j d{\bf x}, \qquad i,j,k=1,\ldots,(n_1\cdots n_d),
 \label{eq:triple}
\end{equation}
and summing them up with the values of ${\bf f}$,
\begin{equation}
\mathcal{L}_{\mathcal{B}}(i,j) = \sum_{k=1}^{n_1\cdots n_d} \mathcal{H}(i,j,k) {\bf f}(k).
\label{eq:bilin-triple}
\end{equation}
Notice that we assume the basis functions can be enumerated by $d$ independent indices, i.e. $i$ is equivalent to $(i_1,\ldots,i_L)$ through \eqref{eq:ten-vec}, and similarly for $j$ and $k$.
The triple elements \eqref{eq:triple} therefore admit a TT decomposition (or even a single Kronecker-product term) similar to \eqref{eq:ttm}.
Now, if the ${\bf f}$ tensor can also be approximated in the TT format \eqref{eq:tt}, the bilinear form \eqref{eq:bilin-abs}
can be represented in this format, with the TT ranks proportional (or equal) to those of ${\bf f}$.
Moreover, the sum in \eqref{eq:bilin-triple} factorizes into individual sums over $k_1,\ldots,k_L$, which can be implemented efficiently block by block.

It remains to compute a TT approximation of ${\bf f}$.
From the previous Newton iteration we are given the TT representation \eqref{eq:btt} for ${\bf z},{\bf c},{\bf p},{\bf q}$.
Hence we can rapidly evaluate any element of the solution components, and afterwards the corresponding value of $f$.
In order to construct a TT approximation to ${\bf f}$ using only a few evaluations of $f$, we use the TT-Cross algorithm \cite{ot-ttcross-2010}.
This is similar to the Alternating Linear Scheme outlined above, except that at each step it draws $r_{m-1} r_m$ fibers of the tensor values in the $m$-th direction in order to populate the $m$-th TT block and prepare the optimized fibers for the next step.
In total it evaluates $\mathcal{O}(Lr^2)$ elements of the tensor, which is feasible under our assumption of small TT ranks.
More robust and rank adaptive generalizations of this algorithm have followed \cite{mo-rectmaxvol-2018,sav-qott-2014,so-dmrgi-2011proc}.

\begin{Rem}
\label{rem:indicator}
Forming the diagonal of the indicator matrix $G_{\lambda}$ in \eqref{GN_ControlConstraints} seems also to be a task for the TT-Cross algorithm.
However, it is likely to perform poorly in this setting, for two reasons.
Firstly, if the discontinuity in a function, e.g. $\max\{0,u-u_+\}$, is not aligned to coordinate axes, the corresponding TT approximation requires very large TT ranks.
This can be seen already in a two-dimensional case: a triangular matrix with all ones in one of the triangles is full-rank.
Secondly, the TT-Cross algorithm is likely to overlook the part of the active set which is not covered by the initial (e.g. random) set of samples.
In order to adapt the sampling fibers, the cross methods require a low discrepancy between adjacent tensor elements, which is not the case for $G_{\lambda}$.
For this reason, we apply the TT approach only to the case of the unconstrained control.
\end{Rem}

\section{Numerical experiments with the low-rank approximations}\label{sec:NumEx2}

In this section, we benchmark the TT algorithm and compare it to the solver with the full vector representation.
The initial distribution of bacterial cells $z_0$ and the desired state $\widehat{z}$ are chosen as in \eqref{eq:z0_full} and \eqref{eq:zhat}, with initial and target concentrations for the chemoattractant $c$ set to zero. In this section the model is solved with an unconstrained control $u$, and final time $T=1$.
For the TT computations we used the TT-Toolbox implementation (see https://github.com/oseledets/TT-Toolbox).

\subsection{Benchmarking of full and low-rank solvers}
\begin{figure}[t]
\caption{CPU time (sec.) (left) and TT ranks (right) for different grid sizes $n$, $m_0=3$ initial peaks, accuracy threshold $\eps=10^{-4}$.}
\label{fig:ttimes_n}
\resizebox{!}{0.49\linewidth}{%
\begin{tikzpicture}
\begin{axis}[%
ylabel=CPU time, axis y line=left, y label style={at={(-0.05,1.0)},anchor=south west,rotate=-90}, ymode=log, yminorgrids,%
xlabel=$\log_{2}n$,%
legend style={at={(0.01,0.99)},anchor=north west}]

\addplot+[blue,solid] coordinates{(5,19.0398)
                                  (6,148.73)
                                  (7,1459.5)
                                 }; \addlegendentry{full};
\addplot+[blue,dashed,domain=5:9,no marks] {2^(3.3*x-12.5)}; \addlegendentry{$n^{3.3}$};
\addplot+[red,solid] coordinates{(5,11.498  +     39.297)
                                 (6,20.595  +      74.48)
                                 (7,48.761  +     187.68)
                                 (8,112.28  +     644.47)
                                 (9,250.85  +     2964.1)
                                }; \addlegendentry{TT};
\addplot+[red,dashed,domain=5:9,no marks] {2^(1.7*x-3.8)}; \addlegendentry{$n^{1.7}$};
\end{axis}
\end{tikzpicture}
}\resizebox{!}{0.49\linewidth}{%
\begin{tikzpicture}
\begin{axis}[%
ylabel=TT rank, axis y line=left, y label style={at={(-0.05,1.0)},anchor=south west,rotate=-90},%
xlabel=$\log_{2}n$,%
legend style={at={(0.01,0.99)},anchor=north west}]

\addplot+[blue,solid] coordinates{(5,56.429)
                                  (6,67)
                                  (7,77.571)
                                  (8,83.286)
                                  (9,94.857)
                                };
\end{axis}
\end{tikzpicture}
}
%
%
%
\end{figure}
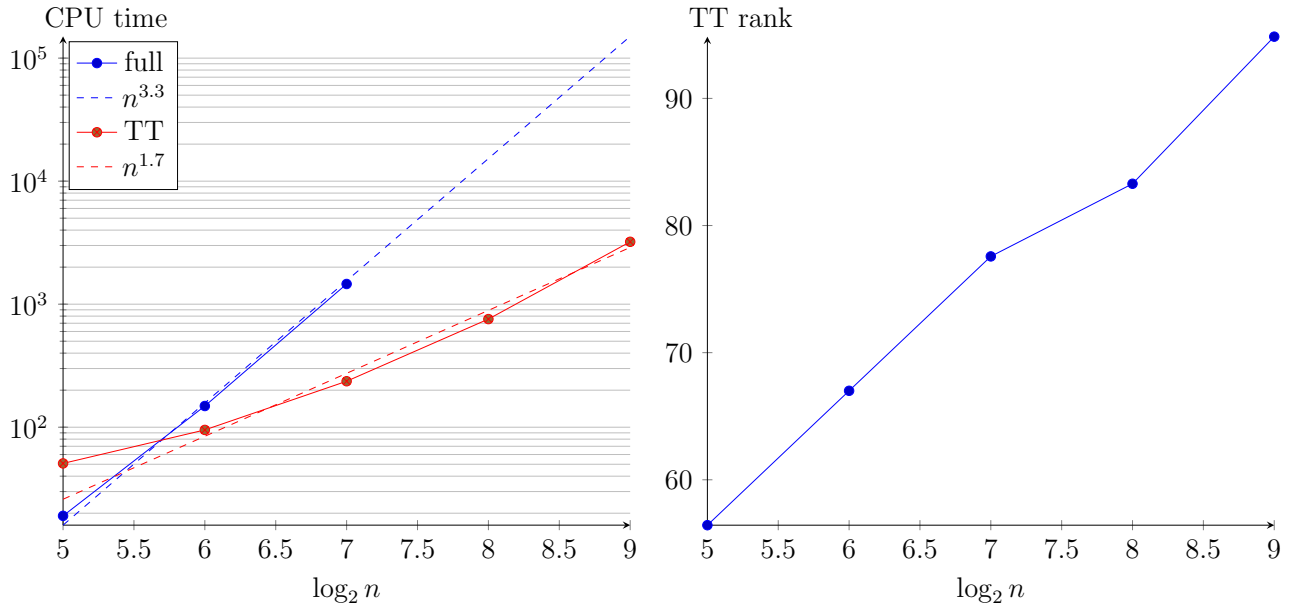

First, we compare CPU times of the original scheme that stores full vectors with those of the approximate TT solver (see Fig. \ref{fig:ttimes_n}).
We fix $m_0=3$ randomly positioned Gaussian peaks in the initial distribution $z_0$.
Since the particular ranks and numbers of iterations depend on the choice of $z_0$,
we average the results over $8$ realizations of $z_0$, for each value of $n$.

The cost of the full-format solver grows slightly faster than cubically, which is expected for a three-dimensional problem.
This concerns both the CPU time and the memory.
In particular, we could not run the full solver for $n>128$ due to the memory limitations.
On the other hand, the TT solver can proceed to much finer grids with lower time and memory footprint.

\subsection{Discretization and TT approximation errors}
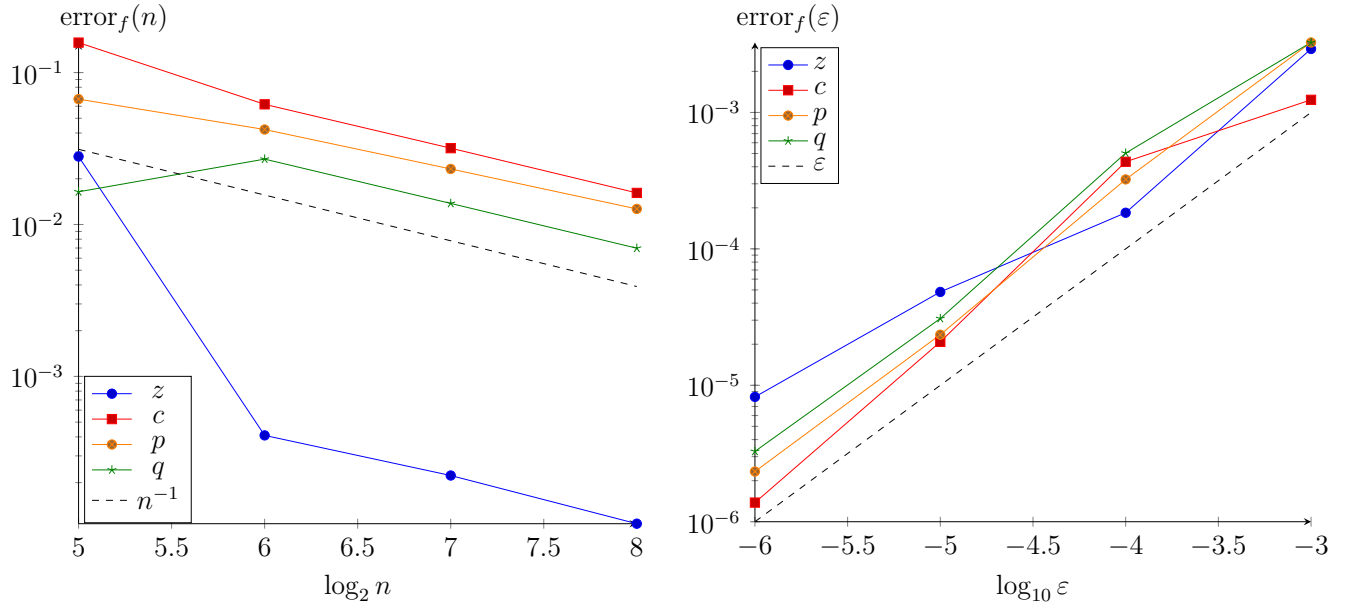
\begin{figure}[t]
\caption{Left: discretization errors for different grid sizes $n$ and $\eps=10^{-6}$. Right: approximation errors for different thresholds $\eps$ and $n=64$.}
\label{fig:err_n_eps}
\resizebox{!}{0.49\linewidth}{
\begin{tikzpicture}
\begin{axis}[%
ylabel=$\mbox{error}_f(n)$, axis y line=left, y label style={at={(-0.05,1.0)},anchor=south west,rotate=-90}, ymode=log, %
xlabel=$\log_{2}n$,%
legend style={at={(0.01,0.001)},anchor=south west}]


\addplot+[blue,solid] coordinates{(5,2.8012e-02)
                                  (6,4.0916e-04)
                                  (7,2.2276e-04)
                                  (8,1.0731e-04)
                                }; \addlegendentry{$z$};
\addplot+[red, solid] coordinates{(5,1.5744e-01)
                                  (6,6.1785e-02)
                                  (7,3.1790e-02)
                                  (8,1.6139e-02)
                                }; \addlegendentry{$c$};
\addplot+[orange,solid] coordinates{(5,6.6834e-02)
                                    (6,4.2184e-02)
                                    (7,2.3201e-02)
                                    (8,1.2664e-02)
                                }; \addlegendentry{$p$};
\addplot+[green!50!black,solid] coordinates{(5,1.6413e-02)
                                            (6,2.6934e-02)
                                            (7,1.3763e-02)
                                            (8,6.9795e-03)
                                }; \addlegendentry{$q$};
\addplot+[black,dashed,domain=5:8,no marks] {2^(-x)}; \addlegendentry{$n^{-1}$};
\end{axis}
\end{tikzpicture}
}
\resizebox{!}{0.49\linewidth}{
\begin{tikzpicture}
\begin{axis}[%
ylabel=$\mbox{error}_f(\eps)$, axis y line=left, y label style={at={(-0.05,1.0)},anchor=south west,rotate=-90}, ymode=log,%
xlabel=$\log_{10}\eps$,%
legend style={at={(0.01,0.99)},anchor=north west}]

\addplot+[blue,solid] coordinates{(-3,2.9210e-03)
                                  (-4,1.8400e-04)
                                  (-5,4.8340e-05)
                                  (-6,8.2210e-06)
                                }; \addlegendentry{$z$};
\addplot+[red, solid] coordinates{(-3,1.2380e-03)
                                  (-4,4.3490e-04)
                                  (-5,2.0820e-05)
                                  (-6,1.3810e-06)
                                }; \addlegendentry{$c$};
\addplot+[orange,solid] coordinates{(-3,3.2570e-03)
                                    (-4,3.2300e-04)
                                    (-5,2.3500e-05)
                                    (-6,2.3330e-06)
                                }; \addlegendentry{$p$};
\addplot+[green!50!black,solid] coordinates{(-3,3.2600e-03)
                                            (-4,5.0360e-04)
                                            (-5,3.0860e-05)
                                            (-6,3.2780e-06)
                                }; \addlegendentry{$q$};
\addplot+[black,dashed,domain=-6:-3,no marks] {10^x}; \addlegendentry{$\eps$};

\end{axis}
\end{tikzpicture}
}
\end{figure}

In order to justify the use of very fine grids (up to $n=512$), let us estimate the discretization errors.
In Fig. \ref{fig:err_n_eps} (left), we vary the grid levels and plot the relative difference of the solutions on the grids with $n$ and $2n$ points in each direction,
$$
\mbox{error}_f(n) = \frac{\left|\langle f^2_{n}\rangle - \langle f_{2n}^2\rangle\right|}{\langle f_{2n}^2\rangle}, \qquad \langle f^2_n \rangle = f_n(T)^\top M f_n(T),
$$
where $f_n(T)$ is the final-time snapshot of the solution component $f \in \{z,c,p,q\}$, computed at the grid with $n$ nodes in each variable,
and $M$ is the mass matrix in space.
The number of initial peaks $m_0=3$ and their positions are fixed in these experiments.

%

We see that the error decays linearly with respect to $n$, as expected from the implicit Euler scheme, for all quantities.
Since this decay is rather slow, at least $256$ points in each direction are necessary to achieve an accuracy of $1\%$ in $q$, and hence the control, $u$.

The truncated singular value decomposition in the TT algorithm tries to introduce the same average amount of error to all solution components.
However, the relative error in each component may differ from $\eps$, depending on the norm scale and other factors of the algorithm, such as the local system solver.
In Fig. \ref{fig:err_n_eps} (right) we investigate the relative error in all components $f\in\{z,c,p,q\}$,
$$
\mbox{error}_f(\eps) = \frac{\|f_{\eps}-f_{10^{-8}}\|_F}{\|f_{10^{-8}}\|_F},
$$
where $f_{\eps}$ is the solution vector computed with the TT approximation threshold $\eps$.
We see that, on average, the errors decay linearly with $\eps$, as expected.



\subsection{Number of peaks in the initial distribution}
\begin{figure}[t]
\caption{CPU time (sec.) (left) and TT ranks (right) for different numbers $m_0$ of initial peaks, accuracy threshold $\eps=10^{-4}$, $n=256$.}
\label{fig:ttimes_m0}
\resizebox{!}{0.49\linewidth}{%
\begin{tikzpicture}
\begin{axis}[%
ylabel=CPU time, axis y line=left, y label style={at={(-0.05,1.0)},anchor=south west,rotate=-90}, ymode=log, ymin=6e2,ymax=2e4,%
xlabel=$m_0$, xmin=3, xmax=20,  %
legend style={at={(0.01,0.99)},anchor=north west}, %
stack plots=y,area style,enlarge x limits=false]
\addplot+[stack plots=false,black,dashed,domain=3:20,no marks] {2^(3.3*8-12.5)};
\node[anchor=south west] at (axis cs: 3,1.5e4) {extrapolation of the full solver time};

\addplot+[fill=white,draw=black] coordinates{
                                     (3, 756.7416-98.6747)
                                     (5, 1.2772e+03-167.1532)
                                     (7, 1.7725e+03-166.4400)
                                     (10,3.1910e+03-411.9074)
                                     (15,7.0464e+03-1.4976e+03)
                                     (20,1.3313e+04-2.0444e+03)
                                } \closedcycle;
\addplot+[draw=blue,fill=red!30!white,mark=*,mark options={fill=blue}] coordinates{
                                  (3, 98.6747)
                                  (5, 167.1532)
                                  (7, 166.4400)
                                  (10,411.9074)
                                  (15,1.4976e+03)
                                  (20,2.0444e+03)
                                } \closedcycle;
\addplot+[draw=black,fill=red!30!white] coordinates{
                                  (3, 98.6747)
                                  (5, 167.1532)
                                  (7, 166.4400)
                                  (10,411.9074)
                                  (15,1.4976e+03)
                                  (20,2.0444e+03)
                                } \closedcycle;
\end{axis}
\end{tikzpicture}
}
\resizebox{!}{0.49\linewidth}{%
\begin{tikzpicture}
\begin{axis}[%
ylabel=TT rank, axis y line=left, y label style={at={(-0.05,1.0)},anchor=south west,rotate=-90},%
xlabel=$m_0$,%
legend style={at={(0.01,0.99)},anchor=north west}, %
stack plots=y,area style,enlarge x limits=false]
\addplot+[fill=white,draw=black] coordinates{
                                     (3, 83.2857-5.8228)
                                     (5, 98-2.2913)
                                     (7, 104.4444-4.7726)
                                     (10,113.2500-4.7434)
                                     (15,121.8889-3.2956)
                                     (20,129.1111-3.6893)
                                } \closedcycle;
\addplot+[draw=blue,fill=red!30!white,mark=*,mark options={fill=blue}] coordinates{
                                  (3, 5.8228)
                                  (5, 2.2913)
                                  (7, 4.7726)
                                  (10,4.7434)
                                  (15,3.2956)
                                  (20,3.6893)
                                } \closedcycle;
\addplot+[draw=black,fill=red!30!white] coordinates{
                                  (3, 5.8228)
                                  (5, 2.2913)
                                  (7, 4.7726)
                                  (10,4.7434)
                                  (15,3.2956)
                                  (20,3.6893)
                                } \closedcycle;
\end{axis}
\end{tikzpicture}
}
\end{figure}
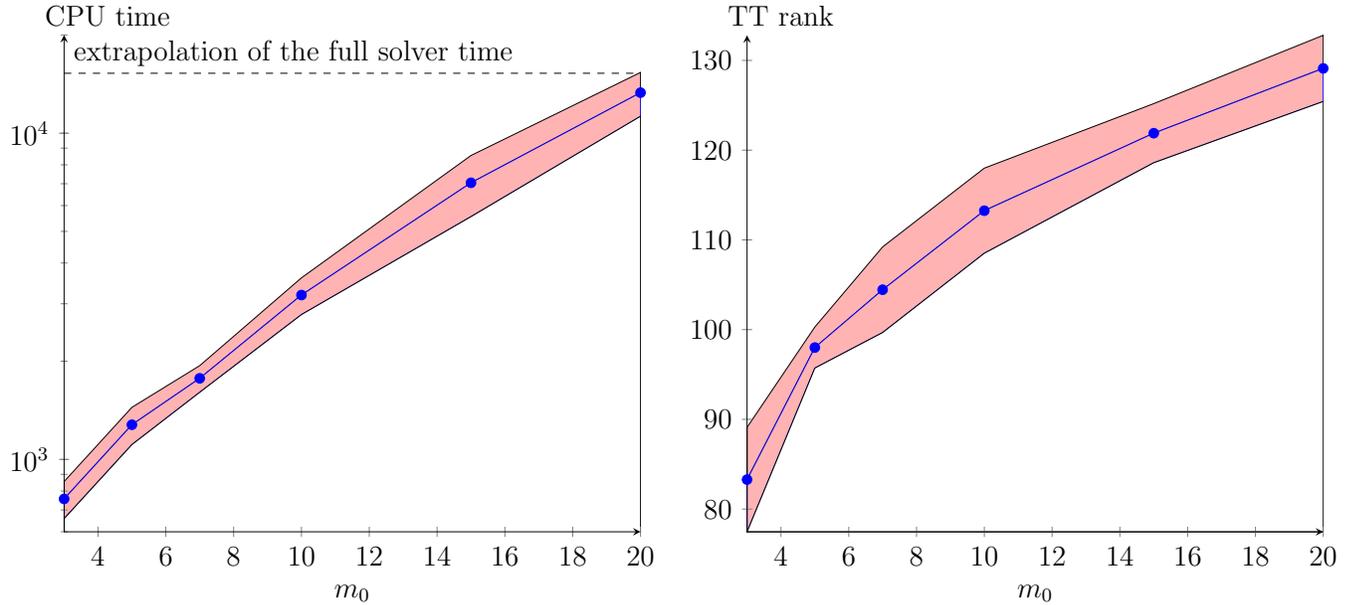

Since the initial distribution of cells \eqref{eq:z0_full} consists of several randomly located Gaussian peaks,
the particular positions of the peaks may influence the performance of the methods.
In Fig. \ref{fig:ttimes_m0} we investigate CPU times and TT ranks in the TT solver versus the number of peaks $m_0$
and their positions.
The plots show means plus minus standard deviations of the times and ranks with respect to the randomization of peak locations.

As expected, the complexity grows with the number of peaks, and for $m_0 \sim 20$ this approaches the estimated time of the full solver (should one have a sufficient amount of memory to run the latter).
For a smaller number of peaks the TT solver is more efficient.
Moreover, a small relative dispersion shows that it is quite insensitive to the particular realization of the initial distribution.

\begin{figure}[htb]
 \centering
 \caption{Left: initial cell density $z_0({\bf x})$ for $m_0=10$. Right: control $u({\bf x},T/2)$.}
 \label{fig:lr_U}
 \includegraphics[width=0.49\linewidth]{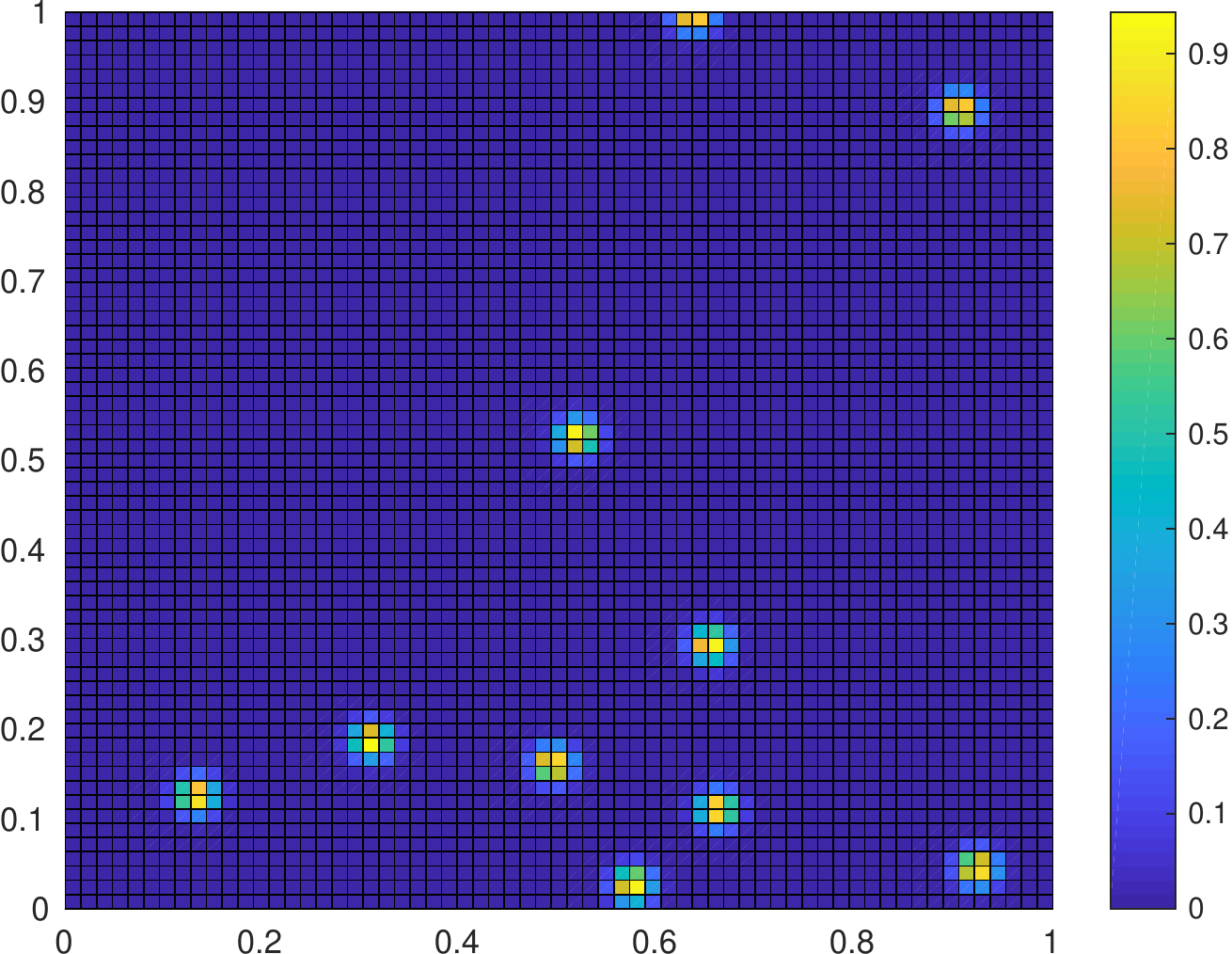}
 \includegraphics[width=0.49\linewidth]{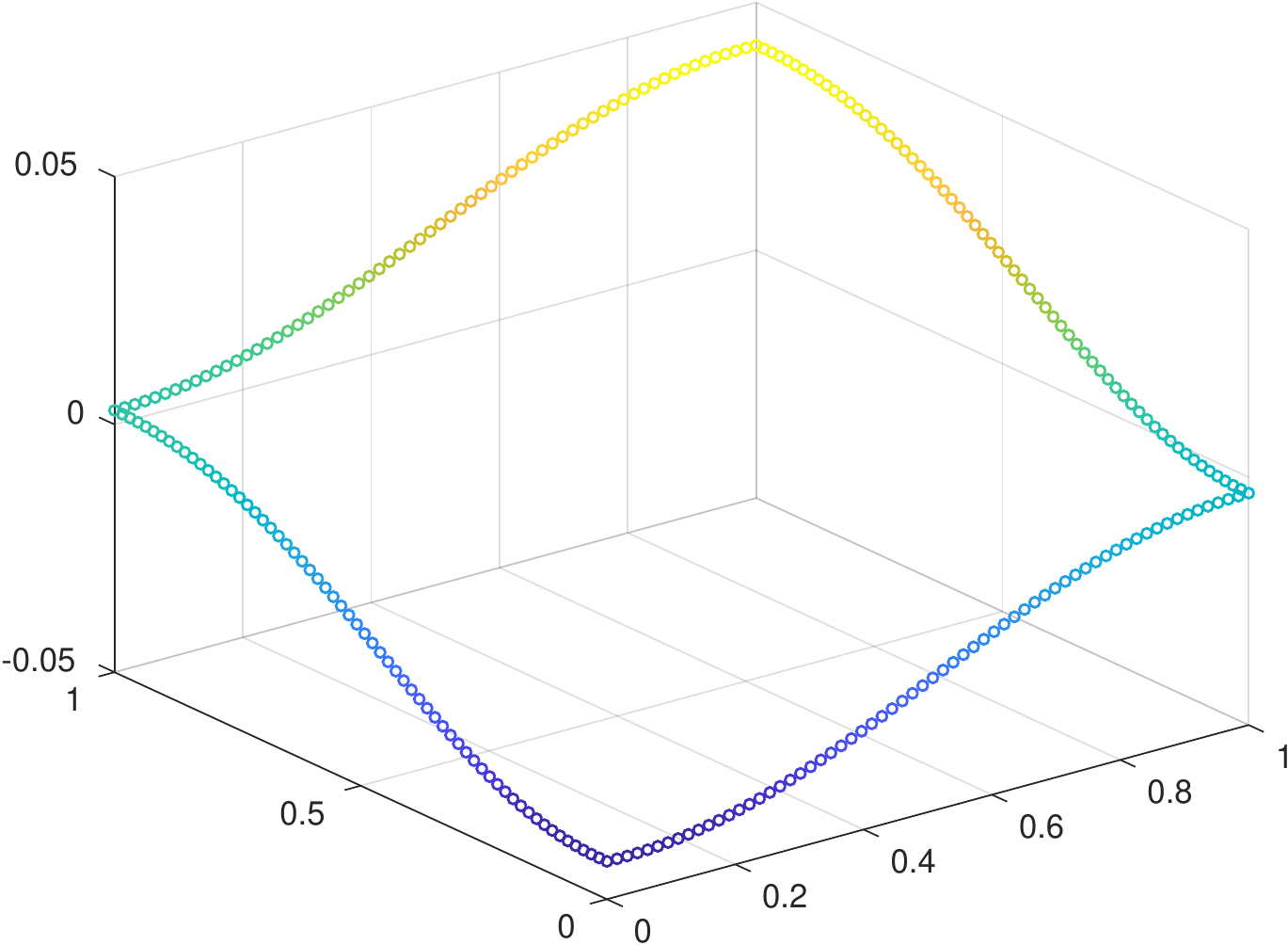}
\end{figure}

\begin{figure}[htb]
 \centering
\caption{Left: cell density at the final time $z({\bf x},T)$. Right: misfit of the density, $z({\bf x},T) - \widehat z({\bf x})$.}
 \label{fig:lr_Z}
 \includegraphics[width=0.49\linewidth]{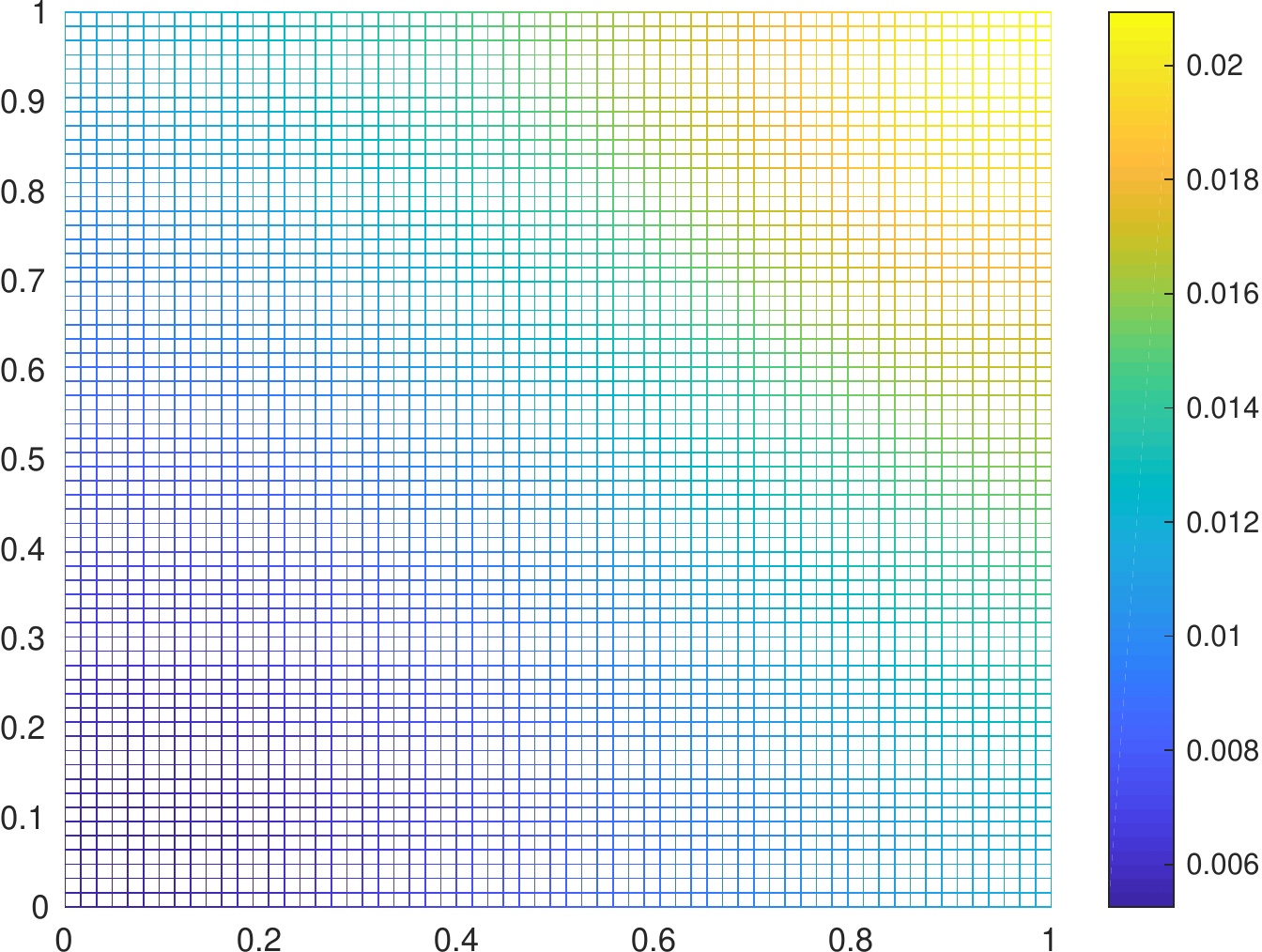}
 \includegraphics[width=0.49\linewidth]{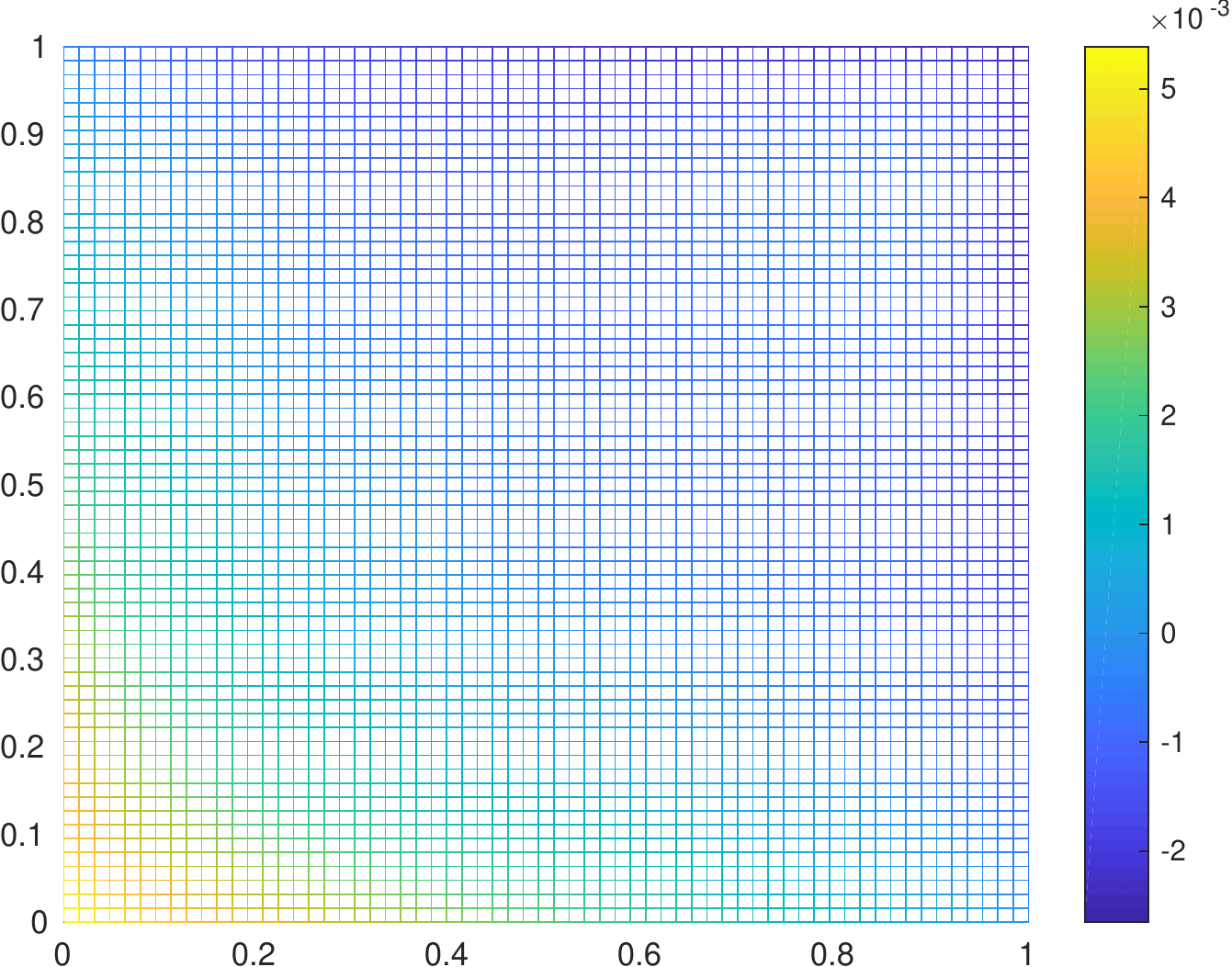}
\end{figure}

The initial cell density for $m_0=10$ peaks and the transient control signal are shown in Fig. \ref{fig:lr_U} (left and right, respectively),
while the final density and the misfit are shown in Fig. \ref{fig:lr_Z}.
The unconstrained control takes negative values in the left bottom corner of the domain.
However, this gives a more accurate fit of the cell density to the desired distribution than the constrained control.

\section{Concluding remarks}\label{sec:Conc}

We have developed a preconditioned Gauss--Newton method for solving optimal control problems in chemotaxis, making use of an effective saddle point type preconditioner coupled with a suitable approximation of the pivoted Schur complement. This enables us to solve potentially huge-scale matrix systems, both without and with additional box constraints imposed on the control variable. Numerical results indicate considerable robustness with respect to the matrix dimension, as well as the parameters involved in the problem set-up.



Moreover, we have shown that the problem without box constraints is amenable to a faster solution using
the low-rank tensor approximations of all vectors and matrices arising in the discretization.
The nonlinearity of the problem can easily be tackled via cross approximation methods, provided that the functions are smooth.
The low-rank decompositions are not very suitable for discontinuous functions, such as an indicator of an active set,
arising in the problem of finding a constrained control or state.
However, in the unconstrained case the low-rank algorithms are much faster and need much less memory than the straightforward solution of the Gauss--Newton equations.
Depending on the ``complexity'' of the transient solution (and hence the tensor ranks), we can achieve a speedup of more than an order of magnitude.

The importance of the box constraints depends on the particular model.
For example, if we can only control the inflow of the chemoattractant, it is reasonable to request a nonnegative control.
However, if the laboratory setup allows one also to remove the chemoattractant, or to add a repellent,
the negative control becomes physically realizable.
This can provide a better control of the cell population,
whereas the low-rank numerical algorithms allow a fast simulation of the required profile of the attractant/repellent, even on a low performance desktop.

\vspace{1em}
\textbf{Acknowledgements.}~~SD and JWP gratefully acknowledge support from the Engineering and Physical Sciences Research Council (UK) Fellowships EP/M019004/1 and EP/M018857/2, respectively.

\bibliographystyle{plain}
\bibliography{refs}

\end{document}